\documentclass[10pt,a4paper]{amsart}

\usepackage[utf8]{inputenc}
\usepackage[T1]{fontenc}
\usepackage[english]{babel}
\usepackage{amsmath, hyperref}
\usepackage[all, 2cell]{xy}
\UseAllTwocells

\title{Simplicial Structure on Complexes}
\author{Djalal Mirmohades}

\newtheorem{theorem}{$\quad$Theorem}[section]
\newtheorem{proposition}[theorem]{$\quad$Proposition}
\newtheorem{corollary}[theorem]{$\quad$Corollary}

\newtheorem{conjecture}[theorem]{$\quad$Conjecture}
\newtheorem{remark}[theorem]{$\quad$Remark}
\newtheorem{example}[theorem]{$\quad$Example}
\theoremstyle{definition}
\newtheorem{definition}[theorem]{$\quad$Definition}
\let\Oldenddefinition\enddefinition
\def\enddefinition{\hfill $\bigtriangleup$\Oldenddefinition}

\renewcommand{\S}{{\mathcal S}}
\newcommand{\A}{{\mathcal A}}

\newcommand{\U}{{\mathcal U}}
\newcommand{\C}{{\mathrm C}}

\newcommand{\HH}{{\mathrm H}}
\newcommand{\N}{{\mathbb N}}
\newcommand{\Q}{{\mathbb Q}}
\newcommand{\Z}{{\mathbb Z}}
\newcommand{\G}{{\mathrm G}}
\newcommand{\F}{{\mathrm F}}
\newcommand{\BG}{\bar{\G}}
\newcommand{\BF}{\bar{\F}}
\newcommand{\li}{{\mathrm{L}}}
\newcommand{\iid}{{\mathrm{\textit{id}}}}
\newcommand{\id}{\mathrm{id}}
\newcommand{\iidd}{{\mathrm{\bf id}}}
\newcommand{\Id}{\mathrm{Id}}
\newcommand{\Cat}{\mathrm{\bf Cat}}
\newcommand{\Set}{\mathrm{\bf Set}}
\newcommand{\Grp}{\mathrm{\bf Grp}}
\newcommand{\Grpd}{\mathrm{\bf Grpd}}
\newcommand{\Sesq}{\mathrm{\bf Sesq}}
\newcommand{\LN}{\operatorname{\bf LN}}
\newcommand{\NN}{\operatorname{\bf N}}
\newcommand{\Ho}{\operatorname{\bf K}}
\newcommand{\ob}{\operatorname{Ob}}
\newcommand{\dom}{\operatorname{\bf dom}}
\newcommand{\cod}{\operatorname{\bf cod}}
\newcommand{\Hom}{\operatorname{Hom}}
\newcommand{\Fun}{\operatorname{Fun}}
\newcommand{\LFun}{\operatorname{LFun}}
\newcommand{\LArr}{\operatorname{LArr}}
\newcommand{\Spine}{\operatorname{Spine}}
\newcommand{\Filler}{\operatorname{Filler}}
\newcommand{\Cone}{\operatorname{Cone}}

\renewcommand{\.}{\kern 0.04em}
\newcommand{\op}[1]{{#1}^{\mathrm{op}}}
\newcommand{\Com}[1]{\mathrm{Com}_{#1}}
\newcommand{\CCom}[1]{\overline{\mathrm{Com}}_{#1}}
\newcommand{\rvec}[2]{\left( \begin{smallmatrix} #1 & #2 \end{smallmatrix} \right)}
\newcommand{\cvec}[2]{\left( \begin{smallmatrix} \vphantom{d} #1 \\ \vphantom{d} #2 \end{smallmatrix} \right)}
\newcommand{\mtrx}[4]{\left( \begin{smallmatrix} \vphantom{\hat{f}} #1 & #2 \\ \vphantom{\hat{f}} #3 & #4 \end{smallmatrix} \right)}
\newcommand{\ord}[1]{{[#1]}}
\newcommand{\nsection}[1]{\subsection{#1} $\:$ \par}

\begin{document}
\begin{abstract}
While chain complexes are equipped with a differential $d$ satisfying $d^2 = 0$, their generalizations called $N$-complexes have a differential $d$ satisfying $d^N = 0$.
In this paper we show that the lax nerve of the category of chain complexes is pointwise adjoint equivalent to the d{\'e}calage of the simplicial category of $N$-complexes.
This reveals additional simplicial structure on the lax nerve of the category of chain complexes which provides a categorfication of the triangulated homotopy category of chain complexes.
We study this phenomena in general and present evidence that the axioms of triangulated categories have simplicial origin.
\end{abstract}
\maketitle

\thispagestyle{empty}
\tableofcontents

\section{Introduction}

\nsection{Categorification}
Decategorification is a special process by which one turns an $n$-category into an $(n-1)$-category.
The reverse process, that is to construct something that decategorifies to a given category, is called \textit{categorification}.
To some extent, categorification aims to introduce new interesting structure which could reveal
previously hidden information.

One way of decatigorifying an $n$-category, which we denote by $\Ho$, is to first identify all isomorphic $(n-1)$-morphisms and then to forget about the $n$-morphisms.
The following example illustrates this.

Let $\Q\!-\!\mathrm{mod}$ be the category of finite-dimensional vector spaces over the rational numbers. 
Then there is a bijection
$$\dim : \Ho(\Q\!-\!\mathrm{mod}) \longrightarrow \N$$
which sends a vector space $V$ to $\dim(V)$. The structure which is inherited in this case is that given by direct sum and tensor product, that is
\begin{align*}
\dim(X \oplus Y) &= \dim(X) + \dim(Y), \\
\dim(X \otimes Y) &= \dim(X) \dim(Y).
\end{align*}
For more examples, see \cite{Maz}.

\nsection{\texorpdfstring{$N$-complexes: an informal description of the results}{N-complexes: an informal description of the results}}
Let $\Com{N}(\A)$, where $N$ is a positive integer, denote the sesquicategory (i.e.\ a \mbox{2-category} without the interchange law) of $N$-complexes over an additive category $\A$, that is ``chain'' complexes that satisfy
$$d^N = 0$$
instead of the usual $d^2 = 0$. There is also an  appropriate generalization of the notion of homotopy which is used to define 2-morphisms.
The category of ordinary chain complexes over $\A$ is then denoted $\Com{2}(\A)$.
The decategorification $\Ho(\Com{2}(\A))$ of $\Com{2}(\A)$, known as the \textit{homotopy category of chain complexes}, has a very special structure called \textit{triangulation}. But this structure emerges only after decategorification.
One serious inconvenience in  the theory of triangulated categories is failure of functoriality in the axiom which 
addresses extensions of morphisms to triangles, see \cite[p. 245]{GM}.
If we go back to $\Com{2}(\A)$, the corresponding construction called the \textit{mapping cone} is a functor, but the presumed triangles
$$0 \longrightarrow \operatorname{dom}(f) \longrightarrow \operatorname{Cyl}(f) \longrightarrow \Cone(f) \longrightarrow 0$$
where $f$ is a morphism, do not satisfy anything that resembles triangulation.
However, in this paper we will show that some of the axioms of triangulated categories, here called \textit{semitriangulation}, can be categorified in a neat way, in terms of a section of a d{\'e}calage which we call \textit{weak recalage}.
Moreover, for the category $\Ho(\Com{2}(\A))$ we find a corresponding (strict) recalage called \textit{the simplicial sesquicategory of $N$-complexes over $\A$} which is denoted $\Com{\Delta}(\A)$.
It is a 2-functor from the opposite simplex 2-category $\op{\Delta}$ to the (large) category of sesquicategories.
It maps an object $\ord{n}$ of $\op{\Delta}$ to the sesquicategory $\Com{n+1}(\A)$, while 1-morphisms and 2-morphisms of $\op{\Delta}$ are mapped to functors and natural transformations, respectively, in a nerve-like manner (see Definition \ref{delta-complex}).

The existence of an object $\Com{\Delta}(\A)$ with its strict interpretation of the axioms of semitriangulated categories was somewhat unexpected.
We believe that there is a deeper reason as to why this happens and that $\Com{\Delta}(\A)$ deserves to be called a categorification of the triangulated category $\Ho(\Com{2}(\A))$.

\nsection{Outline}
In Chapter \ref{category_theory}, we define category theoretic notions such as sesquicategories, decategorification, 2-simplicial objects, the 2-nerve $\NN$, the lax nerve $\LN$, the d{\'e}calage $\Sigma$ and the 2-simplicial functor
$$\Spine_\S : \S \longrightarrow \NN(\S_{\ord{0}}).$$

In Chapter \ref{n_complexes}, the simplicial sesquicategory of $N$-complexes $\Com{\Delta}(\A)$ is defined and it is proven, by explicit construction, that $\Spine_{\Sigma \Com{\Delta}(\A)}$ is a pointwise adjoint equivalence.

In Chapter \ref{triangulations}, we argue that the above equivalence shows that $\Com{\Delta}(\A)$, with its nerve-like simplicial structure alone, categorifies $\Ho(\Com{2}(\A))$ together with a significant part of its triangulated structure.

In Chapter \ref{groups}, we go in the opposite direction and show that the decategorification of semitriangulation produces a group structure. Using our experience from the previous chapter, we obtain that a group structure on a set is the same thing as a recalage of the nerve of the set.

In Chapter \ref{conclusion} we summarize our results and present some conclusions and reflections.

\nsection{Some background}
This work began as an attempt to explain Theorem 1.3 in \cite{Ka} by Kapranov which says that the homology of a 3-complex is exact.
As is pointed out in \cite{Ka}, homologies of $N$-complexes for various $N$ are functors $\mathrm{\bf H}$
\[\xymatrix{
\hdots \ar[r]^(.4){\mathrm{\bf H}} & \Com{3}(\A) \ar[r]^{\mathrm{\bf H}} & \Com{2}(\A) \ar[r]^{\mathrm{\bf H}} & \Com{1}(\A)
}\]
that satisfy $\mathrm{\bf H} \circ \mathrm{\bf H} = 0$ at the rightmost part of the diagram.

But why does this hold only for 3-complexes?
An examination of a broader range of homologies in \cite{Mi} by the author suggested that the individual homology functors ${}_pH_i$ of \cite{Ka} are indeed the right ones.
The functor $\mathrm{\bf H}$ is defined as the total complex of the double complex made up of all the ${}_pH_i$, which has nontrivial sums only when $N \geq 4$.
We conclude that addition should be avoided.
To answer the question, we look behind homology into homotopy, for a property that holds for all $N$-complexes, but at the same time explains the above special feature of 3-complexes.
This paper proposes the following version of an answer:

We show, using an explicit recursive construction, that the simplicial functor $\Spine_{\Sigma \Com{\Delta}(\A)}$ is a degreewize adjoint equivalence.
This equivalence can be used to move the forgotten simplicial functors of $\Sigma \Com{\Delta}(\A)$ to a weak recalage of the lax nerve of $\A$.
In particular at degree 2, the equivalence says that 3-complexes are triangles and their homology is hence exact.

The fact that 3-complexes are triangles was mentioned by Alexei Bondal to my advisor Volodymyr Mazorchuk (we did not manage to find any reference for this property, in this paper this appears as Corollary \ref{triangles_threecomplex}).
On the other hand, the homology of the above ``mapping cone''-recalage of the nerve of $\A$ was described by \cite[p. 98]{Wal} in 1966. 
More recently, a derived equivalence was shown in \cite[Th. 4.9]{IKM} which is related to our results; the derived equivalence follows from the adjoint equivalence mentioned above. 
This is because some of the homology functors ${}_pH_i$ factors trough the functors $\mathrm{\bf v}_i$ (these are the components of $\Spine$ see Definition \ref{vertex} and Definition \ref{spine}) and the remaining homology functors are determined by those. Hence the quasi-isomorphisms of $\Sigma \Com{\Delta}(\A)_\ord{n}$ and $\LN_+(\Com{2}(\A))_\ord{n}$ coincide and the adjoint equivalence in Theorem \ref{adj_eq2} induces an equivalence between the derived categories.

The subject of $N$-complexes can be traced back to 1942 in \cite{Maye} by Mayer.
Interest started to grow after the preprint \cite{Ka} by M. Kapranov which appeared in 1991, see for example 
\cite{1,2,3,4,5,6,7,8,9,10,11,12,13} and references therein. Interest seems now to be accelerating due to the paper \cite{KQ} by Khovanov and Qi.

\section{Category Theory}\label{category_theory}

\nsection{Categories}
Terminology is mainly due to \cite{Mac} and \cite{NL}. 

\begin{definition}\label{category}
A \textit{category} consists of a set $O$ of \textit{objects}, a set $A$ of \textit{arrows} or \textit{morphisms}, three functions
\[\xymatrix@C=25pt{A \ar@<1.2ex>[r]^(.45){\mathrm{dom}} \ar@<-1.2ex>[r]_(.45){\mathrm{cod}} & O \ar[l]|{\id}}\]
and for any pair of arrows $f, g$ such that $\mathrm{dom}(f) = \mathrm{cod}(g)$, there exist a unique arrow denoted $fg$, called their \textit{composition}, satisfying the following axioms:
\begin{align*}
\mathrm{dom}(\id(x)) &= x, \\
\mathrm{cod}(\id(x)) &= x, \\
\mathrm{dom}(gf) &= \mathrm{dom}(f), \\
\mathrm{cod}(gf) &= \mathrm{cod}(g), \\
(fg)h &= f(gh), \\
\id(x) f &= f, \\
f \. \id(x) &= f
\end{align*}
for all arrows $f, g, h$ and objects $x$ where the given composition is defined.
\end{definition}

\begin{definition}
A \textit{functor} $F$ is a structure-preserving map between categories. More precisely, it is two functions (both denoted by $F$):
\[\xymatrix@C=25pt{
A \ar[d]_F \ar@<1.2ex>[r]^(.45){\mathrm{dom}} \ar@<-1.2ex>[r]_(.45){\mathrm{cod}} & O \ar[l]|{\id} \ar[d]^F \\
A' \ar@<1.2ex>[r]^(.45){\mathrm{dom}} \ar@<-1.2ex>[r]_(.45){\mathrm{cod}} & O' \ar[l]|{\id}
} \]
that commute with $\mathrm{dom}$, $\mathrm{cod}$, $\id$ and preserve composition;
$$F(fg) = F(f) F(g)$$
for composable arrows $f, g \in A$.
\end{definition}

\begin{example}\label{preorder}
Any preorder $(P, \leq)$, that is a set $P$ with a reflexive and transitive binary relation $\leq \: \subset P \times P$, can be regarded as a category 
\[\xymatrix@C=25pt{\leq \ar@<1.2ex>[r]^(.45){\mathrm{dom}} \ar@<-1.2ex>[r]_(.45){\mathrm{cod}} & P \ar[l]|{\id}}\]
where $\mathrm{dom}$ and $\mathrm{cod}$ are the two projections $\leq \: \subset P \times P \to P$ and $\mathrm{id} : P \to \: \leq$ is the diagonal map.
In this category there is precisely one morphism $x \to y$ if $x \leq y$ and none otherwise.
The reflexivity of the relation provides identity morphisms and transitivity provides composition of morphisms.
\end{example}

To avoid set-theoretical problems we fix a \textit{Grothendieck universe} $\U$ (see \cite{Wi}) and adopt a terminology relative to $\U$.

\begin{definition}
A set $M$ is said to be \textit{small} if $M \in \U$.
\end{definition}

\begin{definition}
Define $\Set$ as the category whose objects are all small sets and morphisms are functions between sets.
\end{definition}

\begin{definition}
A \textit{2-category} is a category $C$ with additional data:
\begin{itemize}
 \item For each pair of objects $(x, y)$ in $C$, there is a category denoted $C(x, y)$ with its object set equal to the set of morphisms $x \to y$. 
 \item For each triple $(x, y, z)$ of objects in $C$, there is a functor called \textit{horizontal composition}
 $$\_ \circ_0 \_ : C(y, z) \times C(x, y) \longrightarrow C(x, z)$$
 that coincides with composition in $C$ on objects $(f, g)$ of $C(y, z) \times C(x, y)$ i.e.\ $f \circ_0 g = fg : x \to z$.
\end{itemize}
 A morphisms of the category $C(x, y)$ is called a \textit{2-morphisms} of the 2-category $C$.
 
\end{definition}

\begin{definition}
Define $\Cat$ as the 2-category whose objects are all small categories, whose morphisms are functors between small categories and whose 2-morphisms are natural transformations between those.
\end{definition}

\begin{definition}\label{object}
There is a functor
$$\ob : \Cat \longrightarrow \Set$$
that maps each small category $C$ to the set $\ob(C)$ of its objects and functors between categories to the corresponding functions.
This functor has both a left and a right adjoint, both of which are sections of $\ob$ (they map a set $M$ to a category $C$ such that $\ob(C) = M$).

The left adjoint of $\ob$ maps a set to a category with no morphisms but the identity morphisms and maps functions to functors in the obvious way.
We will use the left adjoint of $\ob$ implicitly when a set is regarded as a category or when an $n$-category is regarded as an $(n + 1)$-category.

The right adjoint of $\ob$ maps a set to the preorder given by the full binary relation on the set, considered as a category as in Example \ref{preorder}, and maps functions to functors in the obvious way.
\end{definition}

\begin{definition}\label{sesquicategory}
A \textit{sesquicategory} (or 1.5-category) is a category $C$ together with a lift of the hom-functor along $\ob : \Cat \to \Set$. That is a functor $\mathrm{HOM}_C : \op{C} \times C \to \Cat$ such that the diagram
\[\xymatrix{
& \Cat \ar[d]^\ob \\
\op{C} \times C \ar[r]_(.55){\Hom_C} \ar[ru]^(.55){\mathrm{HOM}_C} & \Set
}\]
commutes.
A functor between sesquicategories is an ordinary functor $F : C \to D$ together with a lift of the hom-natural transformation
\begin{align*}
\Hom_F : \Hom_C & \Longrightarrow \Hom_D \circ \: (\op{F} \times F) \\
(\Hom_F)_{(x, y)} : \Hom_C (x, y) & \longrightarrow \Hom_D (Fx, Fy) \\
f & \longmapsto (\Hom_F)_{(x, y)}(f) := Ff \\
\end{align*}
\[\xymatrix@R=20pt{
\op{C} \times C \ar[dr]^{\Hom_C} \ar[dd]_{\op{F} \times F} & \ar@{}[ddl]_(.6){\dir{=>}}|{\:\: \Hom_F} \\
& \Set \\
\op{D} \times D \ar[ur]_{\Hom_D} &
}\]
along $\ob : \Cat \to \Set$. That is a natural transformation
$$\mathrm{HOM}_F : \mathrm{HOM}_C \Longrightarrow \mathrm{HOM}_D \circ (\op{F} \times F)$$
\[\xymatrix@R=20pt{
\op{C} \times C \ar[dr]^{\mathrm{HOM}_C} \ar[dd]_{\op{F} \times F} & \ar@{}[ddl]_(.6){\dir{=>}}|{\:\:\: \mathrm{HOM}_F} \\
& \Cat \ar[r]^\ob & \Set\\
\op{D} \times D \ar[ur]_{\mathrm{HOM}_D}
}\]
such that $\ob \circ_0 \: \mathrm{HOM}_F = \Hom_F$.
A natural transformation between functors is an ordinary natural transformation $\eta : F \Rightarrow G : C \to D$ which is compatible with the chosen lift.
That is, for every pair of objects $x, y \in C$, the diagram
\[\xymatrix@C=65pt{
\mathrm{HOM}_C(x, y) \ar[r]^{\mathrm{HOM}_G(x, y)} \ar[d]_{\mathrm{HOM}_F(x, y)} & \mathrm{HOM}_D(Gx, Gy) \ar[d]^{\mathrm{HOM}_D(\eta_x, \id)} \\
\mathrm{HOM}_D(Fx, Fy) \ar[r]_{\mathrm{HOM}_D(\id, \eta_y)} & \mathrm{HOM}_D(Fx, Gy)
}\]
commutes.
\end{definition}

In a sesquicategory $C$, each hom-set $\Hom(x, y)$ is the object set of a small category $\mathrm{HOM}(x, y)$.
Objects and morphisms of the category $\mathrm{HOM}(x, y)$ are called 1- and 2-morphisms of $C$, respectively.
There is a ``vertical'' composition of 2-morphisms (the composition inside $\mathrm{HOM}(x, y)$) and a ``horizontal'' composition between a 2-morphism $\alpha$ (morphism of $\mathrm{HOM}(x, y)$) and a 1-morphism $f : y \to z$ defined by
$$f \alpha := \mathrm{HOM}(\id_x, f) (\alpha).$$
This composition becomes distributive due to the functoriallity of $\mathrm{HOM}(\id_x, f)$.
That is, given
\vspace{-2em}
\[\xymatrix@C=40pt{
x \ar[r] \rtwocell<8>{\omit} \rtwocell\omit{<-2>^\alpha} \rtwocell\omit{<2>^\beta} & y \ar[r]^f & z \\
}\] $\.$
\vspace{-2em} \\
we have
$$f (\beta \alpha) = (f \beta) (f \alpha).$$
Composition with $f$ from the right is defined analogously.

Note that no ``horizontal'' compositions of 2-morphisms is postulated.
However, one can construct a horizontal composition in two ways; given a diagram
\[\xymatrix@C=40pt{
\bullet \rtwocell<4>^e_f{\alpha} & \bullet \rtwocell<4>^g_h{\beta} & \bullet \\
}\]
we may construct
$$(\beta f) (g \alpha) \mbox{ and } (h \alpha) (\beta e)$$
but unfortunately, these two compositions differ in general.

\begin{proposition}
A 2-category is a sesquicategory where the two horizontal compositions are equal. That is, $(\beta f) (g \alpha) = (h \alpha) (\beta e)$ in the above diagram.
\end{proposition}

The following example is the main reason we need sesquicategories:
\begin{example}
The category of chain complexes over an additive category (see Section \ref{n_complexes}), with chain maps as 1-morphisms and chain homotopies as 2-morphisms, is a sesquicategory (but not a 2-category).
\end{example}

\begin{definition}
Define $\Sesq$ as the 2-category of all small sesquicategories. Morphisms are functors between sesquicategories and 2-morphisms are natural transformations. We will see in Example \ref{ex_ses} that there are also 3-morphisms in this category.
\end{definition}

\begin{definition}\label{sestercategory}
A \textit{sestercategory}\footnote{From Latin s\={e}mis (``half'') and ter (``three times'').} (or 2.5-category) is a 2-category $C$ together with a lift of the hom-functor along $\ob : \Sesq \to \Cat$, that is a 2-functor $\mathrm{HOM}$
\[\xymatrix{
& \Sesq \ar[d]^\ob \\
\op{C} \times C \ar[r]_(.55){\Hom} \ar[ru]^{\mathrm{HOM}} & \Cat
}\]
such that the diagram commutes.
\end{definition}

\begin{example}\label{ex_ses}
The category $\Sesq$ can be regarded as a sestercategory in which the 3-morphisms are modifications.
\end{example}

\begin{definition}
When $C$ is a sesquicategory and $D$ a small category, define
$$\Fun(D, C)$$
as a sesquicategory with objects, morphisms and 2-morphisms given by functors, natural transformations and modifications, respectively.
\end{definition}

\nsection{Equimorphisms}

\newcommand{\eqv}{\tilde{\: \Rightarrow}\:}
\newcommand{\veqv}{\phantom{\sim}\text{\normalsize $\Downarrow$}\sim}

In the following definition, we are interested in at most 2.5-categories, for a definition of 3-categories and higher see \cite{Le}.

\begin{definition}\label{equimorphism}
In an $n$- or $\left( n-\frac{1}{2} \right)$-category, given a pair of parallel $m$-morphisms $x$ and $y$, an \textit{equimorphism} $x \eqv y$ is, for $m = n$
$$\mbox{the equality relation,}$$
for $0 \leq m < n$ a pair of $(m+1)$-morphisms
\[\xymatrix@C=25pt{x \ar@<.6ex>[r]^f & y \ar@<.6ex>[l]^g }\]
together with a pair of equimorphisms
$$fg \eqv \id_y , \: gf \eqv \id_x.$$
All objects in the category are regarded as parallel $0$-morphisms.
For $n$-morphisms $x, y$ we have $x \eqv y \Leftrightarrow x = y$ (Example \ref{preorder} with the equality relation).
\end{definition}

The definition of an equimorphism is recursive. Expanding the recursion, an equimorphism $x \eqv y$ is two $(m+1)$-morphisms
\[\xymatrix@C=40pt{
\bullet \rtwocell<4>^{x}_{y}{\omit} \ar@{}[r]|{f \. \text{\normalsize $\Downarrow \Uparrow$} g} & \bullet \:,
}\]
four $(m+2)$-morphisms 
\[\xymatrix@C=40pt{
y \rtwocell<4>^{fg}_{\id_y}{\omit} \ar@{}[r]|{\epsilon \. \text{\normalsize $\Downarrow \Uparrow$} \zeta} & y &
x \rtwocell<4>^{gf}_{\id_x}{\omit} \ar@{}[r]|{\theta \. \text{\normalsize $\Downarrow \Uparrow$} \eta} & x \:, \\
}\]
eight $(m+3)$-morphisms
\[\xymatrix@C=40pt{
fg \rtwocell<4>^{\zeta \epsilon}_{\id_{fg}}{\omit} \ar@{}[r]|{\text{\normalsize $\Downarrow \Uparrow$}} & fg & \cdots \:,
}\]
and so on up to the $n$-morphisms which must be isomorphisms (see Proposition \ref{isomorphism}).

We say that $x$ and $y$ are \textit{equivalent} if there exist an equimorphism between them.
The equimorphism $x \eqv y$ can be described by the data $(f, g, \epsilon, \eta, \zeta, \theta, \dots)$.
We call $(g, \epsilon, \eta, \zeta, \theta, \dots)$ \textit{a weak inverse of $f$} and we loosely say that $f$ and $g$ are \textit{weak inverses} of each other, implying the existence of a tail $(\epsilon, \eta, \zeta, \theta, \dots)$.

Functors preserve equimorphisms. The functor $\mathrm{HOM}(\id, r)$ (see Definition \ref{sesquicategory}) is used to define a \textit{horizontal composition} between equimorphisms and morphisms
\[\xymatrix@C=40pt{
\bullet \rtwocell<4>^{x}_{y}{\omit} \ar@{}[r]|\veqv & \bullet \ar[r]^r & \bullet
}\]
as it maps an equimorphism $x \eqv y$ to $rx \eqv ry$.
Given equimorphisms $(f, g, \dots) : x \eqv y$ and $(f', g', \dots) : y \eqv z$ between parallel $m$-morphisms $x, y, z$, define their \textit{vertical composition} in the following way:

For $m = n$, the statements $x = y$ and $y = z$ imply
$$x = z.$$

For $m < n$, the pair of $(m+1)$-morphisms are
\[\xymatrix@C=60pt{
\bullet \rtwocell<5>^{x}_{z}{\omit} \ar@{}[r]|{f'f \. \text{\normalsize $\Downarrow \Uparrow$} gg'} & \bullet \:,
}\]
the equimorphism $f'fgg' \eqv \id_z$ is given by applying the functor $\mathrm{HOM}(g', f')$ to the equimorphism $fg \eqv \id_y$ and vertically compose with $f'g' \eqv \id_z$:
\[\xymatrix@C=40pt@R=17pt{
z \xlowertwocell[0,3]{}<-15>_{\id_z}{\omit} \ar[r]^{g'} & y \rtwocell<4>^{fg}_{\id_y}{\omit} \ar@{}[r]|\veqv & y \ar[r]^{f'} & z \\
& \ar@{}[r]|\veqv & &
}\]
\vspace{-3em} \\
and the equimorphism $gg'f'f \eqv \id_x$ is given analogously:
\[\xymatrix@C=40pt@R=17pt{
x \xlowertwocell[0,3]{}<-15>_{\id_x}{\omit} \ar[r]^{f} & y \rtwocell<4>^{g'f'}_{\id_y}{\omit} \ar@{}[r]|\veqv & y \ar[r]^{g} & x \\
& \ar@{}[r]|\veqv & & \\
& & & .
}\]
\vspace{-3em} \\
Note that vertical composition is defined recursively, just like the definition of equimorphism.

\begin{proposition}\label{isomorphism}
In an $n$- or $\left( n-\frac{1}{2} \right)$-category, an equimorphism between two $(n-1)$-morphisms $x$ and $y$, is a pair of $n$-morphisms
\[\xymatrix@C=25pt{x \ar@<.6ex>[r]^f & y \ar@<.6ex>[l]^g }\]
such that $fg = \id_y$ and $gf = \id_x$.
\end{proposition}
Hence for $1$-categories, equimorphisms between objects are just isomorphisms.

\begin{definition}\label{adj_eq}
A tuple $(F, G, \epsilon, \eta)$ of functors
\[\xymatrix@C=25pt{C \ar@<.6ex>[r]^F & D \ar@<.6ex>[l]^G }\]
and natural transformations (or more generally lax natural transformations, see Definition \ref{LFun-ob})
\[\xymatrix@C=40pt{
D \rtwocell<4>^{FG}_{\id_D}{\epsilon} & D &
C \rtwocell<4>^{GF}_{\id_C}{^\eta} & C \:, \\
}\]
in $\Sesq$ is said to be an \textit{adjoint equivalence} if it can be extended to an equimorphism
$$(F, G, \epsilon, \eta, \dots) : C \eqv D$$
and satisfies the ``triangular identities'' (see \cite[IV.1(9)]{Mac})
$$\epsilon_{F} (F \eta) = \id_{\F} , \: (G \epsilon) \eta_{G} = \id_{G}.$$
\end{definition}

\nsection{Categorification and Decategorification}

\begin{definition}
A \textit{groupoid} is a category where all morphisms are isomorphisms.
Define $\Grpd$ as the full subcategory of $\Cat$ which contains all groupoids.

\end{definition}

\begin{definition}
Define a functor
$$\operatorname{iso} : \Cat \longrightarrow \Grpd$$
that maps each category $C$ to the groupoid $\operatorname{iso}(C)$ consisting of the same objects as $C$ but having only the isomorphisms of $C$ as morphisms.
Obviously $\operatorname{iso}(C)$ has identity morphisms and is closed under composition of morphisms.
Functors between categories are mapped to their restriction to the isomorphisms, this works because functors preserve isomorphisms.
\end{definition}

\begin{definition}
Define a functor
$$\pi_0 : \Grpd \longrightarrow \Set$$
that maps each groupoid $G$ to the set of isomorphism classes $\pi_0(G)$ of $G$ and functors between groupoids to functions between sets in the obvious way.
\end{definition}

\begin{definition}\label{def_decat}
\textit{Decategorification} is a functor $\Ho : \Cat \to \Set$ defined as
$$\Ho := \pi_0 \circ \operatorname{iso}.$$
We extend this functor to $\Ho : \Sesq \to \Cat$ by applying it to the underlying morphism categories.
By \textit{double decategorification} we mean the composition $\Ho \circ \Ho : \Sesq \to \Set$.

We also consider decategorification of internal theories (see \cite[XII.1]{Mac} for an example) in $\Cat$ and $\Sesq$.
That is, structures defined in terms of categories, functors, natural transformations (and modifications) that satisfy certain relations.
When doing so, we require that
$$\Ho \mbox{preserves equimorphisms.}$$
This essentially means that whenever a theory requires two things to be isomorphic, then the decategorified theory requires the decategorified things to be equal.
We say that $A$ is a \textit{categorification} of $B$ whenever $B$ is a decategorification of $A$.
\end{definition}

\nsection{Lax Structures}

\begin{definition}
 Given a sesquicategory $C$, define the sesquicategory $\LArr(C)$, called \textit{the lax arrow category of $C$}, together with three functors 
 \[\xymatrix@C=35pt{\LArr(C) \ar@<1.2ex>[r]^(.6){\dom} \ar@{<-}[r]|(.6){\iidd} \ar@<-1.2ex>[r]_(.6){\cod} &C}\]
 as follows:
 
 \textbf{An object} $x$ of $\LArr(C)$ is a 1-morphism of $C$. Define $\dom x := \mathrm{dom}(x)$ and $\cod x := \mathrm{cod}(x)$; the domain and codomain of $x$ in $C$.
 
  \textbf{A 1-morphism} $f : x \to y$ in $\LArr(C)$ is a triple $(f_0, f_1, \widehat{f})$ consisting of:
  \begin{itemize}
   \item A 1-morphism $f_0 : \mathrm{dom}(x) \to \mathrm{dom}(y)$ in $C$.
   \item A 1-morphism $f_1 : \mathrm{cod}(x) \to \mathrm{cod}(y)$ in $C$.
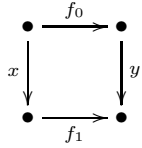
 \[\xymatrix{
  \bullet \ar[d]_x \ar[r]^{f_0} & \bullet \ar[d]^y \\
  \bullet \ar[r]_{f_1} & \bullet \\
 }\]
   \item A 2-morphism $\widehat{f} : y f_0 \Rightarrow f_1 x$ in $C$.
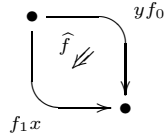
 \[\xymatrix{
  \bullet \ar`r[r][dr]^(-.5){y f_0} \ar`^r[d][dr]_(-.5){f_1 x} & \ar@{}[dl]_{\widehat{f}}|{\dir{=>}} \\
  & \bullet \\
 } \]
  \end{itemize}
 Define $\dom f := f_0$ and $\cod f := f_1$.
 
  \textbf{A 2-morphism} $\alpha : f \to g : x \to y$ in $\LArr(C)$ is a pair $(\widehat{\alpha}_0, \widehat{\alpha}_1)$ consisting of:
  \begin{itemize}
   \item A 2-morphism $\widehat{\alpha}_0 : f_0 \to g_0$ in $C$.
   \item A 2-morphism $\widehat{\alpha}_1 : f_1 \to g_1$ in $C$.
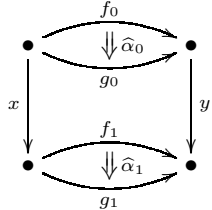
 \[\xymatrix@R=35pt@C=50pt{
  \bullet \ar[d]_x \rtwocell^{f_0}_{g_0}{\:\:\: \widehat{\alpha}_0} & \bullet \ar[d]^y \\
  \bullet \rtwocell^{f_1}_{g_1}{\:\:\: \widehat{\alpha}_1} & \bullet \\
 }\]
  \end{itemize}
 that satisfy the following equality
$$\widehat{g} \circ_1 (y \circ_0 \widehat{\alpha}_0) = (\widehat{\alpha}_1 \circ_0 x) \circ_1 \widehat{f} .$$
 Define $\dom \alpha := \widehat{\alpha}_0$ and $\cod \alpha := \widehat{\alpha}_1$.

 For an object $x$ in $\LArr(C)$, the identity morphism $\id_x$ of $\LArr(C)$ is the triple
 $(\id_{\mathrm{dom}(x)}, \id_{\mathrm{cod}(x)}, \iid_x)$:
 \[\xymatrix{
  \bullet \ar[d]_x \ar[r]^{\id_{\mathrm{dom}(x)}} & \bullet \ar[d]^x & & \bullet \ar`r[r]^(1.4){x}[dr] \ar`^r[d][dr]_(-.4){x} & \ar@{}[dl]_{\iid_{x}}|{\dir{=>}} \\
  \bullet \ar[r]_{\id_{\mathrm{cod}(x)}} & \bullet & & & \bullet \\
 }\]
 where $\iid_x$ is the identity 2-morphism of $x$.
 
The functor $\iidd$ maps an object $a$ of $C$ to the object $\id_{a}$ of $\LArr(C)$, a morphism $f : a \to b$ to the triple $(f, f, \iid_{f})$:
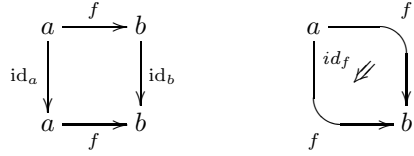
 \[\xymatrix{
  a \ar[d]_{\id_{a}} \ar[r]^{f} & b \ar[d]^{\id_{b}} & & a \ar`r[r]^(1.4){f}[dr] \ar`^r[d][dr]_(-.4){f} & \ar@{}[dl]_{\iid_{f}}|{\dir{=>}} \\
  a \ar[r]_{f} & b & & & b \\
 }\]
and a 2-morphism $\alpha$ to the pair $(\alpha, \alpha)$.
Note the difference between identity morphisms in $\LArr(C)$ and the image of the $\iidd$ functor; the former are made of ``horizontal'' identities while the later ``vertical'' ones.

Composition of 1-morphisms $f$ and $g$ in $\LArr(C)$ is defined as
$$f g := (f_0 g_0, \. f_1 g_1, \. (f_1 \circ_0 \widehat{g})  \circ_1  (\widehat{f} \circ_0 g_0))$$
 
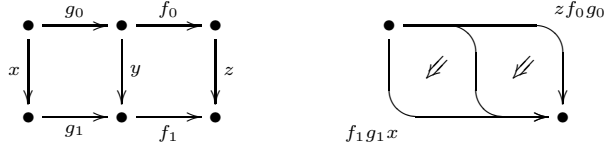
\[\xymatrix{
  \bullet \ar[d]_x \ar[r]^{g_0} & \bullet \ar[d]^y \ar[r]^{f_0} & \bullet \ar[d]^z &&
  \bullet \ar`r[r]`^r[dr][drr] \ar`r[rr]^(1.3){z f_0 g_0}[drr] \ar`^r[d][drr]_(-.3){f_1 g_1 x}
  & \ar@{}[dl]|{\dir{=>}} & \ar@{}[dl]|{\dir{=>}} \\
  \bullet \ar[r]_{g_1} & \bullet \ar[r]_{f_1} & \bullet &&
  & & \bullet\\
 }\]
 where $\circ_0$ and $\circ_1$ denote horizontal and vertical composition in $C$, respectively.
Vertical composition of 2-morphisms is defined pointwise and horizontal composition between a 1-morphisms $f$ and a 2-morphism $\alpha$ is defined as
$$f \circ_0 \alpha := (f_0 \circ_0 \widehat{\alpha}_0, f_1 \circ_0 \widehat{\alpha}_1)$$
and similarly when the order is reversed.
It is straightforward to verify that these definitions satisfy the required equality of 2-morphisms and that vertical composition with 1-morphisms behaves like a functor as required by Definition \ref{sesquicategory}.
\end{definition}

\begin{definition}\label{LFun-ob}
Let $C$ be a sesquicategory and $D$ a small category. Define $\LFun(D, C)$ as the sesquicategory given by:

\textbf{Objects} are all functors $F : D \to C$.
 
\textbf{1-morphisms} $\alpha : F \to G : D \to C$, called \textit{lax natural transformations} are defined as functors
$$\alpha : D \longrightarrow \LArr(C)$$
such that $\dom \circ \: \alpha = F$ and $\cod \circ \: \alpha = G$. 
 
\textbf{2-morphisms} $\eta : \alpha \to \beta : F \to G : D \to C$ , called \textit{modifications} are defined as natural transformations 
$$\eta : \alpha \to \beta : D \longrightarrow \LArr(C)$$
such that $\dom \eta = \id_F$ and $\cod \eta = \id_G$.

For a functor $F : D \to C$, the corresponding identity lax natural transformation $\id_F$ is given by 
$$\id_F := \iidd \circ F : D \longrightarrow \LArr(C).$$

Given two lax natural transformations
\[\xymatrix{
F \ar[r]^\alpha & G \ar[r]^\gamma & H \\
}\]
define their composition $\gamma \alpha$ on objects $x$ of $D$ as
$$(\gamma \alpha) x := (\gamma x) (\alpha x)$$
and on morphisms $f : x \to y$ of $D$
$$(\gamma \alpha) f := ((\alpha f)_0, \: (\gamma f)_1, \: (\iid_{\alpha x} \circ_0 \widehat{\gamma f})  \circ_1  (\widehat{\alpha f} \circ_0 \iid_{\gamma y} ))$$ 
\[\xymatrix{
Fx \ar[d]_{Ff = (\alpha f)_0} \ar[r]^{\alpha x} & Gx \ar[d]|{(\alpha f)_1 = Gf = (\gamma f)_0} \ar[r]^{\gamma x} & Hx \ar[d]^{(\gamma f)_1 = Hf} &&
Fx \ar`r[r]`^r[dr][drr] \ar`r[rr]^(1.5){(Hf) (\gamma x) (\alpha x)}[drr] \ar`^r[d][drr]_(-.5){(\gamma y) (\alpha y) Ff} & & \\
Fy \ar[r]_{\alpha y} & Gy \ar[r]_{\gamma y} & Hy &&
\ar@{}[ur]|(.55){\dir{=>}} & \ar@{}[ur]|(.55){\dir{=>}} & Hy \\
}\]
 
Given a modification $\eta : \alpha \to \beta$ and a lax natural transformation $\gamma$
\[\xymatrix@C=40pt{
F \rtwocell<4>^{\alpha}_{\beta}{\eta} & G \ar[r]^\gamma & H
}\]
we define their composition $\gamma \eta$ on objects $x$ of $D$ as
$$(\gamma \eta)_x := (\gamma x) \eta_x.$$
\end{definition}
 
\begin{definition}\label{LFun-mor}
Let $G: A \to B$ be a functor between small categories and $C$ a sesquicategory. Define the 2-functor
$$\LFun(G, C) : \LFun(B, C) \longrightarrow \LFun(A, C)$$
by composition from the right:

\textbf{Objects} $F$ of $\LFun(B, C)$ (functors $B \to C$) are mapped to $F \circ G$.
 
\textbf{1-morphisms} $\alpha$ of $\LFun(B, C)$ (lax natural transformations) are viewed as functors $\alpha : B \to \LArr(C)$ and are mapped to $\alpha \circ G : A \to \LArr(C)$.
 
\textbf{2-morphisms} $\eta : \alpha \to \beta$ of $\LFun(B, C)$ (modifications) are viewed as natural transformations $\eta : \alpha \to \beta : B \to \LArr(C)$ and are mapped to $\eta_G : \alpha \circ G \to \beta \circ G$.
 
\end{definition}

\begin{proposition}\label{fun_emb}
There is a faithful functor injective on objects
$$\Fun(D, C) \longrightarrow \LFun(D, C)$$
where a functor is mapped to the same functor and a natural transformation $\alpha :$ $F \to G :$ $D \to C$ is mapped to the functor $\alpha' : D \to \LArr(C)$ given by 
$$X \longmapsto \alpha_X : FX \longrightarrow GX$$
for objects $X$ of $D$ and 
$$f \longmapsto (Ff, Gf, \iid) : \alpha_X \longrightarrow \alpha_Y$$
for morphisms $f : X \to Y$, where $\iid$ denotes the identity 2-morphism on $\alpha_Y (Ff) = (Gf) \alpha_X$. Then $\dom \circ \, \alpha' = F$ and $\cod \circ \, \alpha' = G$.
\end{proposition}

\nsection{Simplicial Categories}
Let $\N = \{0, 1, 2, \dots \}$ be the set of natural numbers.

\begin{definition}\label{ordcat}
For $n \in \N \cup \{-1\}$, define $\ord{n}$ as the category given by the preorder (see Example \ref{preorder})
$$\left( \{0, 1, \cdots, n\}, \. \leq \right)$$
where $\leq$ is the order relation between natural numbers.
\end{definition}
In the category $\ord{n}$, there is precisely one arrow $k \to m$ if $k \leq m$, none otherwise.
The category has $n+1$ objects and ${(n+2)(n+1)}/{2}$ morphisms, in particular $\ord{-1}$ is the empty category.

\begin{definition}\label{def_olik}
Define the \textit{simplex 2-category}, denoted $\Delta$, as the full subcategory of $\Cat$ with objects $\{ \ord{n} \, | \, n \in \N \}$. 1-morphisms (functors) are then order-preserving functions $f : \ord{n} \to \ord{m}$.
There exists a unique 2-morphism (natural transformation) denoted $\prec \; : f \to g : \ord{n} \to \ord{m}$ if and only if $\forall x \in \ord{n} : f(x) \leq g(x)$.

\end{definition}

\begin{definition}\label{def_aug}
Define the \textit{augmented simplex 2-category}, denoted $\Delta_+$, as the full subcategory of $\Cat$ with objects $\{ \ord{n} \, | \, n \in \N \cup \{-1\} \}$.
\end{definition}

\begin{definition}
Given a 2-category $C$, define the \textit{opposite of $C$}, denoted $\op{C}$, as the category with the same objects, its 1-morphisms reversed (i.e.\ each 1-morphism in $C$ has its domain and codomain swapped) and its 2-morphisms unchanged:
\[\xymatrix@C=30pt{
x \rtwocell^f_g{\alpha} & y & \mbox{ is replaced with } & x & y \ltwocell^g_f{^\alpha} \: .\\
}\]
\end{definition}

\begin{definition}
A \textit{simplicial category} is a 2-functor $\op{\Delta} \to \Cat$.
\end{definition}

\begin{definition}
An \textit{augmented simplicial category} is a 2-functor $\op{\Delta}_+ \to \Cat$.

\end{definition}

From now on, let $S$ be a fixed simplicial category and use the following shorthand notation.
Let
$$S_{\ord{n}} := S(\ord{n}).$$
Let
$$\operatorname{\bf d}_i := S(d^i : \ord{n-1} \to \ord{n}) : S_{\ord{n}} \longrightarrow S_{\ord{n-1}}, \: \: 0 \leq i \leq n,$$
for the face functors and
$$\mathrm{\bf s}_i := S(s^i : \ord{n+1} \to \ord{n}) : S_{\ord{n}} \longrightarrow S_{\ord{n+1}}, \: \: 0 \leq i \leq n,$$
for the degeneracy functors. Given the 2-morphism
$$\prec \: : f \longrightarrow g : \ord{n} \longrightarrow \ord{m}$$
let the natural transformation
$$\tau = \tau_{f, g} \: := S(\prec) : S(f) \longrightarrow S(g) : S_{\ord{m}} \to S_{\ord{n}}.$$
Objects of $S_{\ord{n}}$ are called $n$-simplicies.

\begin{definition}\label{decalage}
Consider the 2-functor $\sigma : \Delta \to \Delta$ ($\sigma : \Delta_+ \to \Delta_+$) which maps the object $\ord{n}$ to $\ord{n+1}$, morphisms $d^n$ to $d^{n+1}$, $s^{n}$ to $s^{n+1}$ and 2-morphisms in the only possible way.
Now, given an (augmented) simplicial category $S$, define \textit{the d{\'e}calage of $S$}, or $\Sigma S$ as
$$\Sigma S := S \circ \sigma.$$
\end{definition}
One may think of $\Sigma S$ as an (augmented) simplicial category, shifted in degree 
$$(\Sigma S)_{\ord{n}} = S_{\ord{n+1}}$$
where $S_{\ord{0}}$ ($S_{\ord{-1}}$) and all the $\mathrm{\bf d}_0$ and $\mathrm{\bf s}_0$ have been forgotten.
Note that by remembering $S_{\ord{0}}$, $\Sigma S$ becomes augmented, even if $S$ wasn't.

\nsection{Nerve Constructions}

\begin{definition}\label{nerve}
Given a small category $C$, define the simplicial category $\NN(C)$, called the \textit{2-nerve of $C$}, as
$$\NN(C) := \Fun( \: \_ \: , \: C) |_{\op{\Delta}} : \op{\Delta} \to \Cat$$
where $\Fun$ denotes the functor category (internal hom in $\Cat$) and $|_{\op{\Delta}}$ denotes the restriction to $\op{\Delta}$.
\end{definition}

\begin{definition}\label{augnerve}
Analogously, define the simplicial category $\NN_+(C)$, called the \textit{augmented 2-nerve of $C$}, as
$$\NN_+(C) := \Fun( \: \_ \: , \: C) |_{\op{\Delta}_+} : \op{\Delta}_+ \to \Cat.$$
\end{definition}

By the 2-Yoneda lemma we have
$$\mathrm{2 \mbox{\textendash} Nat}(\NN(\ord{n}), \NN(C)) \simeq \NN(C)_{\ord{n}} = \Fun(\ord{n}, \: C).$$

\begin{definition}
Given a simplicial category $S$, define the simplicial set $\ob(S)$, where at degree $n \in \N$
$$\ob(S)_n := \ob(S_{\ord{n}})$$
with the face and degeneracy functions given by the corresponding functors.
\end{definition}

\begin{example}
The simplicial set $\ob(\NN(C))$ is known as \textit{the nerve of $C$} (see \cite[XII.2]{Mac}).
\end{example}

\begin{definition}
Given a small 2-category $C$, define the simplicial category $\big(\LN_+(C)\big)$ $\LN(C)$, called the \textit{lax-nerve of $C$}, as
\begin{align*}
\LN(C) &:= \LFun( \: \_ \: , \: C) |_{\op{\Delta}} : \op{\Delta} \to \Cat \\
\Big( \LN_+(C) &:= \LFun( \: \_ \: , \: C) |_{\op{\Delta}_+} : \op{\Delta}_+ \to \Cat \Big) 
\end{align*}
where $\LFun$ is defined in Definition \ref{LFun-ob} and Definition \ref{LFun-mor}, and $|_{\op{\Delta}}$ denotes the restriction to $\op{\Delta}$.
\end{definition}

\nsection{Simplicial Functors}

\begin{definition}
A \textit{simplicial functor} between two simplicial categories is defined as a natural transformation between their underlying functors.
\end{definition}

\begin{proposition}\label{inc}
There is a faithful simplicial functor injective on objects
$$\NN(C) \longrightarrow \LN(C),$$
induced by the embedding $\Fun \to \LFun$ (see Proposition \ref{fun_emb}).
\end{proposition}

\begin{definition}\label{vertex}
Let $\S$ be a simplicial category and $n$ a positive integer. For each $0 \leq i \leq n$, define the \textit{$i$:th vertex functor}
$$\mathrm{\bf v}_i : \S_{\ord{n}} \longrightarrow \S_{\ord{0}},$$
as the composition
$$\mathrm{\bf v}_i := \mathrm{\bf d}_0 \cdots \mathrm{\bf d}_{i-1} \mathrm{\bf d}_{i+1} \cdots \mathrm{\bf d}_n$$
where $\mathrm{\bf d}_i$ is omitted. The dependence on $n$ is implicitly given by the domain of $\mathrm{\bf v}_i$. For $n = 0$ the functor $\mathrm{\bf v}_0 : \S_{\ord{0}} \to \S_{\ord{0}}$ is the identity functor.
\end{definition}

\begin{proposition}
There are natural transformations 
$$\tau = \S(\prec) : \mathrm{\bf v}_i \to \mathrm{\bf v}_{i+1}.$$
\end{proposition}
The 2-morphism $\prec$ is defined in Definition \ref{def_olik}.
\begin{proof}
We have
$$\prec \: : d_{i+1} \longrightarrow d_i : \ord{i} \longrightarrow \ord{i+1}$$
which is mapped by $\S$ to
$$\tau : \mathrm{\bf d}_{i+1} \longrightarrow \mathrm{\bf d}_i : \S_{\ord{i+1}} \longrightarrow \S_{\ord{i}}.$$
Composing with $\mathrm{\bf d}_0 \cdots \mathrm{\bf d}_{i-1}$ from the left and $\mathrm{\bf d}_{i+2} \cdots \mathrm{\bf d}_n$ from the right we get
$$\tau : \mathrm{\bf d}_0 \cdots \mathrm{\bf d}_{i-1} \mathrm{\bf d}_{i+1} \mathrm{\bf d}_{i+2} \cdots \mathrm{\bf d}_n \longrightarrow \mathrm{\bf d}_0 \cdots \mathrm{\bf d}_{i-1} \mathrm{\bf d}_i \mathrm{\bf d}_{i+2} \cdots \mathrm{\bf d}_n.$$
\end{proof}

\begin{definition}\label{spine}
Given a simplicial category, we shall define a simplicial functor
$$\Spine = \Spine_\S : \S \longrightarrow \NN(\S_{\ord{0}}).$$
Regard the category $\NN(\S_{\ord{0}})_{\ord{n}} = \Fun(\ord{n}, \S_{\ord{0}})$ as the category of diagrams of the shape
$$\xymatrix{0 \ar[r] & 1 \ar[r] & \dots \ar[r] & n}$$
in $\S_{\ord{0}}$ and define the functors $\Spine_{\ord{n}} : \S_{\ord{n}} \longrightarrow \NN(\S_{\ord{0}})_{\ord{n}}$
$$\Spine_{\ord{n}} := \xymatrix{\mathrm{\bf v}_0 \ar[r]^\tau & \mathrm{\bf v}_1 \ar[r]^\tau & \dots \ar[r]^\tau & \mathrm{\bf v}_n}$$
for all $n \geq 0$. For $n = 0$, this means $\Spine_{\ord{0}} = \mathrm{\bf v}_0 : \S_{\ord{0}} \to \S_{\ord{0}}$ is the identity functor.
\end{definition}

\begin{proposition}\label{pspine}
If $\S$ it the 2-nerve of a category, then
$$\Spine : \S \longrightarrow \NN(\S_{\ord{0}})$$
is an isomorphism of simplicial categories.
\end{proposition}

\section{\texorpdfstring{$N$\textendash Complexes}{N-Complexes}}\label{n_complexes}
\nsection{\texorpdfstring{The Simplicial Sesquicategory of $N$\textendash Complexes}{The Simplicial Sesquicategory of N-Complexes}}

\begin{definition}\label{com_n}
For an additive category $\A$ and an integer $N \geq 1$, let $\Com{N}(\A)$ denote the sesquicategory of \mbox{$N$-complexes} over $\A$, defined as follows:

\textbf{An object} $X$ of $\Com{N}(\A)$, called an \textit{$N$-complex} (or \textit{chain complex} when \mbox{$N = 2$}), is a collection of objects $\{X^i\}_{i \in \Z}$ in $\A$ together with a collection of morphisms $d$, called \textit{differentials}, between adjacent pair of objects
\[ \xymatrix@C=35pt@R=35pt{
\cdots \ar[r]^d & X^0 \ar[r]^d & X^1 \ar[r]^d & X^2 \ar[r]^d & \cdots .
} \]
such that the composition of $N$ consecutive differentials equals zero.
The ``degree'' of the differential is never written out, so we may write
$$d^N = 0$$
despite the fact that $d$ is not an endomorphism.

\textbf{A 1-morphism} $f : X \to Y$ of $\Com{N}(\A)$ is a collection of morphisms $\{f^i\}_{i \in \Z}$
\[ \xymatrix@C=35pt@R=35pt{
\cdots \ar[r]^d & X^0 \ar[r]^d \ar[d]^{f^0} & X^1 \ar[r]^d \ar[d]^{f^1} & X^2 \ar[r]^d \ar[d]^{f^2} & \cdots \\
\cdots \ar[r]^d & Y^0 \ar[r]^d & Y^1 \ar[r]^d & Y^2 \ar[r]^d & \cdots
} \]
such that every square commutes (it is essentially a natural transformation).

\textbf{A 2-morphism} $\hat{h} : f \to g : X \to Y$ of $\Com{N}(\A)$
$$\xymatrix{X \rrtwocell<5>^f_g{\hat{h}} & & Y},$$
called a \textit{homotopy}, is a collection of morphisms $\{ \hat{h}^i : X^{i+N-1} \to Y^i\}_{i \in \Z}$ in $\A$ such that
$$g - f = \sum_{i = 1}^N d^{N-i} \hat{h} \. d^{i-1}$$
where the upper indices of the differentials denote exponentiation.
A dotted arrow inside a square
\[\xymatrix@C=30pt@R=30pt{
\bullet \ar[d]_b \ar@{.>}[dr]^(.6){\hat{h}} \xtwocell[dr]{}\omit{_} \ar[r]^a & \bullet \ar[d]^c && \ar@{}[d]|{\mbox{describes the $2$-morphism}} && \bullet \ar`r[r][dr]^(-.3){ca} \ar`^r[d][dr]_(-.3){db} & \ar@{}[dl]_{\hat{h}}|{\dir{=>}} \\
\bullet \ar[r]_d & \bullet &&&& & \bullet}\]

The 1-morphisms behave like natural transformations, composition is defined pointwise and identity 1-morphisms are given by pointwise identities.
Vertical composition of 2-morphisms is written additively and is defined as pointwise addition
$$\hat{h} + \hat{k} = \{ \hat{h}^i + \hat{k}^i \}_{i \in \Z},$$
hence identity 2-morphisms are given by pointwise zero morphisms.
Horizontal composition between 1-morphisms and 2-morphisms is written multiplicatively and is defined as pointwise composition
$$f \hat{h} g = \{ f^i \hat{h}^i g^{i+N-1} \}_{i \in \Z}.$$
\end{definition}

\begin{proposition}\label{com1}
Every object of $\Com{1}(\A)$ is homotopy equivalent to the zero 1-complex. This makes the identity functor on $\Com{1}(\A)$ homotopy equivalent to the (constant) zero functor.
\end{proposition}
\begin{proof}
In $\Com{1}(\A)$, the definitions of morphism and 2-morphism coincide.
\end{proof}

\begin{definition} \label{delta-complex}
Define $\Com{\Delta}(\A)$ to be a simplicial category as follows:

\textbf{On objects} $\ord{n}$ of $\Delta$, define
$$\Com{\Delta}(\A)_{\ord{n}} := \Com{n+1}(\A).$$
An $n$-simplex of $\Com{\Delta}(\A)$ is then an $(n+1)$-complex.

\textbf{On 1-morphisms}, or order-preserving functions $\ord{n} \to \ord{m}$, it is by the lemma in \cite[VII.5]{Mac} enough to define $\Com{\Delta}(\A)$ on the generators $d^i$ and $s^j$.
Given $d^i : \ord{n-1} \to \ord{n}$ in $\Delta$ where $0 \leq i \leq n$, define the functor 
$$\mathrm{\bf d}_i = \Com{\Delta}(\A)(d^i) : \Com{n+1}(\A) \longrightarrow \Com{n}(\A)$$
as the removal of all objects at degree congruent to $i$ modulo $n+1$. A morphisms $f : X \to Y$ in $\Com{n+1}(\A)$ is then mapped to
\[ \xymatrix@R=3pt@C=17pt{
\cdots \ar[r]^d & X^0 \ar[ddd]^{f_0} \ar[r]^d & \cdots \ar[r]^d & X^{i-1} \ar[ddd]^{f_{i-1}} \ar[r]^{d^2} & 
X^{i+1} \ar[ddd]^{f_{i+1}} \ar[r]^d & \cdots \ar[r]^d & X^n \ar[ddd]^{f_n} \ar[r]^d & \cdots \\
& & & & & & & \\
& & & & & & & \\
\cdots \ar[r]^d & Y^0 \ar[r]^d & \cdots \ar[r]^d & Y^{i-1} \ar[r]^{d^2} & Y^{i+1} \ar[r]^d & \cdots \ar[r]^d & Y^n \ar[r]^d & \cdots \\
\cdots & 0 & \cdots & i-1 & i & \cdots & n-1 & \cdots \\
} \]
with the integers at the bottom indicating the new degrees.

Given $s^j : \ord{n+1} \to \ord{n}$, \: $0 \leq j \leq n$, define the functor
$$\mathrm{\bf s}_j = \Com{\Delta}(\A)(s^j) : \Com{n+1}(\A) \longrightarrow \Com{n+2}(\A)$$
as the repetition of objects at degree congruent to $j$ modulo $n+1$ with the identity map in between. A morphisms $f : X \to Y$ in $\Com{n+1}(\A)$ is then mapped to
\[ \xymatrix@R=3pt@C=17pt{
\cdots \ar[r]^d & X^0 \ar[ddd]^{f_0} \ar[r]^d & \cdots \ar[r]^d & X^{j} \ar[ddd]^{f_{j}} \ar[r]^{1} & X^{j} \ar[ddd]^{f_{j}} \ar[r]^d & \cdots \ar[r]^d & X^n \ar[ddd]^{f_n} \ar[r]^d & \cdots \\
& & & & & & & \\ & & & & & & & \\
\cdots \ar[r]^d & Y^0 \ar[r]^d & \cdots \ar[r]^d & Y^{j} \ar[r]^{1} & Y^{j} \ar[r]^d & \cdots \ar[r]^d & Y^n \ar[r]^d & \cdots \\
\cdots & 0 & \cdots & j & j+1 & \cdots & n+1 & \cdots \\
} \]
with the integers at the bottom indicating the new degrees.

\textbf{On 2-morphisms} define the natural transformation
$$\: \tau = \Com{\Delta}(\A)(\prec \: : d^i \to d^{i-1}) : \mathrm{\bf d}_i \longrightarrow \mathrm{\bf d}_{i-1} : \Com{n+1}(\A) \longrightarrow \Com{n}(\A)$$
on objects $X$ of $\Com{n+1}(\A)$, as the morphism $\tau_X : \mathrm{\bf d}_i \: X \longrightarrow \mathrm{\bf d}_{i-1} \: X$ which is the identity map on all degrees, except those degrees congruent to $i-1$ modulo $n$, at which it is the differential $d$ of $X$:
\[ \xymatrix@R=3pt@C=15pt{
&\:\: \cdots \ar[r]^d & X^0 \ar[ddd]^1 \ar[r]^d & \cdots \ar[r]^d & X^{i-2} \ar[ddd]^1 \ar[r]^d & X^{i-1} \ar[ddd]^d \ar[r]^{d^2} & X^{i+1} \ar[ddd]^1 \ar[r]^d & \cdots \ar[r]^d & X^n \ar[ddd]^1 \ar[r]^d & \cdots \\
&\:\: & & & & & & & & \\ & & & & & & & & \\ 
&\:\: \cdots \ar[r]^d & X^0 \ar[r]^d & \cdots \ar[r]^d & X^{i-2} \ar[r]^{d^2} & X^i \ar[r]^d & X^{i+1} \ar[r]^d & \cdots \ar[r]^d & X^n \ar[r]^d & \cdots \\
&\:\: \cdots & 0 & \cdots & i-2 & i-1 & i & \cdots & n-1 & \cdots \\
} \]
Analogously, define
$$\:\:\:\:\:\: \tau = \Com{\Delta}(\A)(\prec \: : s^j \to s^{j+1}) : \mathrm{\bf s}_j \longrightarrow \mathrm{\bf s}_{j+1} : \Com{n+1}(\A) \longrightarrow \Com{n+2}(\A)$$
on objects $X$ of $\Com{n+1}(\A)$, as the morphism $\tau_X : \mathrm{\bf s}_j \: X \longrightarrow \mathrm{\bf s}_{j+1} \: X$
defined as the identity map except at degree congruent to $j+1$ modulo $n+2$, where it is the differential $d$ of $X$:
\[ \xymatrix@R=3pt@C=15pt{
\cdots \ar[r]^d & X^0 \ar[ddd]^1 \ar[r]^d & \cdots \ar[r]^d & X^{j} \ar[ddd]^1 \ar[r]^1 & X^{j} \ar[ddd]^d \ar[r]^{d} & X^{j+1} \ar[ddd]^1 \ar[r]^d & \cdots \ar[r]^d & X^n \ar[ddd]^1 \ar[r]^d & \cdots \\
& & & & & & & & \\ & & & & & & & & \\ 
\cdots \ar[r]^d & X^0 \ar[r]^d & \cdots \ar[r]^d & X^j \ar[r]^d & X^{j+1} \ar[r]^1 & X^{j+1} \ar[r]^d & \cdots \ar[r]^d & X^n \ar[r]^d & \cdots \\
\cdots & 0 & \cdots & \text{\phantom{+}}j\text{\phantom{+}} & j+1 & j+2 & \cdots & n+1 & \cdots \\
} \]
\end{definition}

\begin{definition}\label{com_n_mod}
For an additive category $\A$ and an integer $N \geq 1$, let $\CCom{N}(\A)$ denote the sesquicategory of $\Z/N\Z$-graded $N$-complexes over $\A$.
Equivalently, $\CCom{N}(\A)$ is the subcategory of $\Com{N}(\A)$ which is invariant under degree shift by $\pm N$. That means an object $X$ of $\CCom{N}(\A)$ satisfy
$$X^i = X^j$$
whenever $i \equiv j \pmod{N}$ and accordingly for the differential of $X$.
Similarly for morphisms and 2-morphisms $\CCom{N}(\A)$.
\end{definition}

\begin{definition} \label{delta-complex_mod}
Define the simplicial sesquicategory $\CCom{\Delta}(\A)$ on objects $\ord{n}$ of $\Delta$ as
$$\CCom{\Delta}(\A)_{\ord{n}} := \CCom{n+1}(\A).$$
Morphisms and 2-morphisms are defined analogously to $\Com{\Delta}(\A)$ given in Definition \ref{delta-complex}.
\end{definition}

\begin{theorem}\label{meta_theorem}
Every functor
$$F : \Com{n}(\A) \longrightarrow \Com{m}(\A)\phantom{.}$$
we consider or construct in this paper can be restricted to a functor
$$\overline{F} : \CCom{n}(\A) \longrightarrow \CCom{m}(\A).$$
That means whenever an object $X$ of $\Com{n}(\A)$ happens to be in $\CCom{n}(\A)$, then $FX$ lies in $\CCom{m}(\A)$ and similarly for morphisms and 2-morphisms.
\end{theorem}

\begin{remark}
If $\A$ is not $\U$-small, then the universe $\U$ was not large enough.
In that case, postulate the existence of a bigger universe $\U'$ that contains $\A$ and let $\Cat$ denote the category of $\mathcal{U}'$-small categories.
This makes $\Com{\Delta}(\A)$ a simplicial category.
\end{remark}

\nsection{\texorpdfstring{The Functor $\Filler$}{The Functor Filler}}
\renewcommand{\S}{\widehat{\mathcal A}}

\begin{definition}
Given an additive category $\A$, let
$$\S := \Sigma \Com{\Delta}(\A)$$
or
$$\S := \Sigma \CCom{\Delta}(\A)$$
where $\Com{\Delta}(\A)$ is defined in Definition \ref{delta-complex}, $\CCom{\Delta}(\A)$ is defined in Definition \ref{delta-complex_mod} and $\Sigma$ is defined in Definition \ref{decalage}.
The second case relies on Theorem \ref{meta_theorem}.
\end{definition}

\begin{definition}
Let $\S_{\ord{n}} \, {}_{\mathrm{\bf v}_n} \mkern-8mu \times_{\dom} \LArr(\S_{\ord{0}})$ denote the pullback of the diagram 
\[ \xymatrix@C=20pt@R=20pt{
\S_{\ord{n}} \ar[dr]_{\mathrm{\bf v}_n} & & \LArr(\S_{\ord{0}}) \ar[dl]^{\dom} \\
& \S_{\ord{0}} & \\
} \]
in $\Sesq$. Objects of this category are pairs $(C, x)$, with $C \in \S_{\ord{n}}, x \in \LArr(\S_{\ord{0}})$ that satisfy
$$\mathrm{\bf v}_n C = \dom x.$$
The pair is written $(C, x : \mathrm{\bf v}_n C \to X)$ in order to denote $\cod x$ by $X$.
Although $\mathrm{\bf v}_n C$ is a $2$-complex, its differential will never be written out, instead $d$ will denote the differential of the $(n+2)$-complex $C$ and $\mathrm{\bf v}_n C$ will be written as
\[ \xymatrix@C=35pt@R=35pt{
\cdots \ar[r]^d & C^0 \ar[r]^{d^{n+1}} & C^{n+1} \ar[r]^d & C^{n+2} \ar[r]^{d^{n+1}} & \cdots .
} \]
Because of this convention, every second component of the differential $d$ of $X$ will be decorated with a star
\[ \xymatrix@C=35pt@R=35pt{
\cdots \ar[r]^d & X^0 \ar[r]^{d^\star} & X^1 \ar[r]^d & X^2 \ar[r]^{d^\star} & \cdots .
} \]
The morphism $x : \mathrm{\bf v}_n C \to X$ is a commuting diagram
\[ \xymatrix@C=35pt@R=35pt{
\cdots \ar[r]^d & C^0 \ar[r]^{d^{n+1}} \ar[d]^x & C^{n+1} \ar[r]^d \ar[d]^x & C^{n+2} \ar[r]^{d^{n+1}} \ar[d]^x & \cdots \\
\cdots \ar[r]^d & X^0 \ar[r]^{d^\star} & X^1 \ar[r]^d & X^2 \ar[r]^{d^\star} & \cdots
} \]
and the star is a reminder of the fact that the left square satisfies
$$d^\star x = x d^{n+1}$$
in contrast for the right-hand side square where we have $d x = x d$.
\end{definition}

\begin{definition}
Define the functor
\begin{align*}
\F_n : \S_{\ord{n+1}} & \longrightarrow \S_{\ord{n}} \, {}_{\mathrm{\bf v}_n} \mkern-8mu \times_{\dom} \LArr(\S_{\ord{0}})
\end{align*}
as
$$\F_n := \left( \mathrm{\bf d}_{n+1}, \. \tau \right), \mbox{ where } \tau : \mathrm{\bf v}_n \to \mathrm{\bf v}_{n+1}.$$
To be more precise, $\F_n$ maps an object $C$ of $\S_{\ord{n+1}}$, that is an $n+3$ complex, to the pair consisting of the $n+2$ complex $\mathrm{\bf d}_{n+1} C$
\[ \xymatrix@C=30pt@R=30pt{
\cdots \ar[r]^{d^2} & C^0 \ar[r]^d & C^1\ar[r]^d & \cdots \ar[r]^d & C^{n+1} \ar[r]^{d^2} & C^{n+3}\ar[r]^d & \cdots
} \]
and the morphism of $2$-complexes $\tau_C : \mathrm{\bf v}_n C \to \mathrm{\bf v}_{n+1} C$
\[ \xymatrix@C=35pt@R=35pt{
\cdots \ar[r]^{d^2} & C^0 \ar[r]^{d^{n+1}} \ar[d]^1 & C^{n+1} \ar[r]^{d^2} \ar[d]^d & C^{n+3} \ar[r]^{d^{n+1}} \ar[d]^1 & \cdots \\
\cdots \ar[r]^d & C^0 \ar[r]^{d^{n+2}} & C^{n+2} \ar[r]^d & C^{n+3} \ar[r]^{d^{n+2}} & \cdots .
} \]
regarded as an object of $\LArr(\S_{\ord{0}})$.
Morphisms $f : C \to D$ are mapped to the pair $\left( \mathrm{\bf d}_{n+1} f, \. \tau_f \right)$ where $\tau_f$ is a commuting square in $\S_{\ord{0}}$
\[ \xymatrix@C=35pt@R=35pt{
\mathrm{\bf v}_n C \ar[d]_{\tau_C} \ar[r]^{\tau_{f, 0}} & \mathrm{\bf v}_n D \ar[d]^{\tau_D}\\
\mathrm{\bf v}_{n+1} C \ar[r]_{\tau_{f, 1}} & \mathrm{\bf v}_{n+1} D
} \]
regarded as an morphism of $\LArr(\S_{\ord{0}})$.
\end{definition}

\begin{definition}
Define the functor
\begin{align*}
\G_n : \S_{\ord{n}} \, {}_{\mathrm{\bf v}_n} \mkern-8mu \times_{\dom} \LArr(\S_{\ord{0}}) & \longrightarrow \S_{\ord{n+1}}. 
\end{align*}
as follows:

$\bullet$ On objects $(C, x : \mathrm{\bf v}_n C \to X)$, where $C$ is an $(n + 2)$-complex and $x$ a morphism
\[ \xymatrix@C=35pt@R=35pt{
\cdots \ar[r]^d & C^0 \ar[r]^{d^{n+1}} \ar[d]^x & C^{n+1} \ar[r]^d \ar[d]^x & C^{n+2} \ar[r]^{d^{n+1}} \ar[d]^x & \cdots \\
\cdots \ar[r]^d & X^0 \ar[r]^{d^\star} & X^1 \ar[r]^d & X^2 \ar[r]^{d^\star} & \cdots ,
} \]
define the $(n+3)$-complex $\G_n(C, x)$ as:
\small
\[ \xymatrix@C=11pt@R=20pt{
\cdots \ar[r] & C^0 \oplus X^0 \ar[r]^{\mtrx{d}{0}{0}{1}} & C^1 \oplus X^0 \ar[r]^(.6){\mtrx{d}{0}{0}{1}} & \cdots \ar[r]^(.3){\mtrx{d}{0}{0}{1}} & C^{n+1} \oplus X^0 \ar[r]^(.55){\mtrx{-d}{0}{x}{d^\star}} & C^{n+2} \oplus X^1 \ar[r]^{\mtrx{-1}{0}{x}{d}} & C^{n+2} \oplus X^2 \ar[r]^(.6){\mtrx{d}{0}{0}{1}} & \cdots
} \]
\normalsize
where the diagram starts at degree $0$, ends at degree $n+3$ and the pattern repeats with periodicity $n+3$.

$\bullet$ On morphisms $(a, f) : (C, x) \to (D, y)$
where $f = (f_0, f_1, \hat{f})$, $f_0 = \mathrm{\bf v}_n a$ and $\hat{f} : y f_0 \Rightarrow f_1 x$ is a homotopy
\[\xymatrix@C=35pt@R=35pt{
\mathrm{\bf v}_n C \ar[d]_x \ar[r]^{\mathrm{\bf v}_n a} \ar@{.>}[dr]^(.6){\hat{f}} \xtwocell[dr]{}\omit{_} & \mathrm{\bf v}_n D \ar[d]^y \\
X \ar[r]_{f_1} & Y \\
}\]
\[ \xymatrix@C=35pt@R=35pt{
\cdots \ar[r]^d & C^0 \ar[r]^{d^{n+1}} & C^{n+1} \ar[r]^d \ar[dl]_{\hat{f}} & C^{n+2} \ar[r]^{d^{n+1}} \ar[dl]_{\hat{f}} & \cdots \\
\cdots \ar[r]^d & Y^0 \ar[r]^{d^\star} & Y^1 \ar[r]^d & Y^2 \ar[r]^{d^\star} & \cdots ,
} \]
define the morphism of $(n+3)$-complexes $\G_n(a, f)$ as:
\small
\[ \xymatrix@C=10pt@R=30pt{
\cdots \ar[r] & C^0 \oplus X^0 \ar[r]^{\mtrx{d}{0}{0}{1}} \ar[d]|{\mtrx{a}{0}{\hat{f}d^{n+1}}{f_1}} & C^1 \oplus X^0 \ar[r]^(.65){\mtrx{d}{0}{0}{1}} \ar[d]|{\mtrx{a}{0}{\hat{f}d^n}{f_1}} & \cdots \ar[r]^(.3){\mtrx{d}{0}{0}{1}} & C^{n+1} \oplus X^0 \ar[r]^(.55){\mtrx{-d}{0}{x}{d^\star}} \ar[d]|{\mtrx{a}{0}{\hat{f}}{f_1}} & C^{n+2} \oplus X^1 \ar[r]^{\mtrx{-1}{0}{x}{d}} \ar[d]|{\mtrx{a}{0}{\hat{f}}{f_1}} & C^{n+2} \oplus X^2 \ar[r]^(.67){\mtrx{d}{0}{0}{1}} \ar[d]|{\mtrx{a}{0}{\hat{f}d^{n+1}}{f_1}} & \cdots \\
\cdots \ar[r] & D^0 \oplus Y^0 \ar[r]_{\mtrx{d}{0}{0}{1}} & D^1 \oplus Y^0 \ar[r]_(.65){\mtrx{d}{0}{0}{1}} & \cdots \ar[r]_(.3){\mtrx{d}{0}{0}{1}} & D^{n+1} \oplus Y^0 \ar[r]_(.55){\mtrx{-d}{0}{y}{d^\star}} & D^{n+2} \oplus Y^1 \ar[r]_{\mtrx{-1}{0}{y}{d}} & D^{n+2} \oplus Y^2 \ar[r]_(.67){\mtrx{d}{0}{0}{1}} & \cdots .
} \]
\normalsize

Note that $\G_n$ maps identity morphisms to identity morphisms because in such a case $a = 1$, $f_1 = 1$ and $\hat{f} = 0$. To verify that $\G_n$ preserves composition, recall that composition of morphisms in $\S_{\ord{n}} \, {}_{\mathrm{\bf v}_n} \mkern-8mu \times_{\dom} \LArr(\S_{\ord{0}})$ is defined as 
$$(b, g) \. (a, f) = (b, (g_0, g_1, \hat{g})) \. (a, (f_0, f_1, \hat{f})) = (ba, (g_0 f_0, g_1 f_1, \hat{h}))$$
where $\hat{h} = \hat{g} f_0 + g_1 \hat{f}$ as seen in the diagram
\[\xymatrix@C=35pt@R=35pt{
\mathrm{\bf v}_n C \ar[d]_x \ar[r]^{f_0} \ar@{.>}[dr]^(.6){\hat{f}} \xtwocell[dr]{}\omit{_} & \mathrm{\bf v}_n D \ar[d]^y \ar[r]^{g_0} \ar@{.>}[dr]^(.6){\hat{g}} \xtwocell[dr]{}\omit{_} & \mathrm{\bf v}_n E \ar[d]^z \\
X \ar[r]_{f_1} & Y \ar[r]_{g_1} & Z .\\
}\]
At position $n+2$ modulo $n+3$ we have
\begin{align*}
\G_n(a, f) = \mtrx{a}{0}{\hat{f}}{f_1} &= \mtrx{f_0}{0}{\hat{f}}{f_1}, & \G_n(b, g) = \mtrx{b}{0}{\hat{g}}{g_1} & = \mtrx{g_0}{0}{\hat{g}}{g_1}
\end{align*}
and 
$$\G_n(b, g) \. \G_n(a, f) = \mtrx{g_0}{0}{\hat{g}}{g_1} \mtrx{f_0}{0}{\hat{f}}{f_1} = \mtrx{g_0 f_0}{0}{\hat{g} f_0 + g_1 \hat{f}}{g_1 f_1} = \G_n(ba, gf).$$
At other positions
$$\G_n(b, g) \. \G_n(a, f) = \mtrx{b}{0}{\hat{g} d^k}{g_1} \mtrx{a}{0}{\hat{f} d^k}{f_1} = \mtrx{b a}{0}{\hat{g} d^k a + g_1 \hat{f} d^k}{g_1 f_1} = \G_n(ba, gf)$$
because $\hat{g} d^k a = \hat{g} a d^k = \hat{g} f_0 d^k$.

$\bullet$ On 2-morphisms $(\hat{h}' , \hat{h}'') : (0, 0) \Rightarrow (a, f)$, where $\hat{h}'' = (\hat{h}_0, \hat{h}_1) : (0, 0) \Rightarrow (f_0, f_1)$ is a pair of homotopies
\[\xymatrix@C=60pt@R=50pt{
\bullet \ar[d]_x \rtwocell<3.5>^0_{f_0}{\:\: \hat{h}_0} \ar@{.>}@/_/[dr]^(.6){\hat{f}} \xtwocell[dr]{}\omit{_<1>} & \bullet \ar[d]^y \\
\bullet \rtwocell<3.5>^0_{f_1}{\:\: \hat{h}_1} & \bullet \\
}\]
that satisfy
\begin{align*}
f_0 & = d \hat{h}_0 + \hat{h}_0 d , \\
f_1 & = d \hat{h}_1 + \hat{h}_1 d , \\
\hat{f} & = \hat{h}_1 x - y \hat{h}_0
\end{align*}
such that
$$\hat{h}_0 = \mathrm{\bf v}_n \hat{h}' : \mathrm{\bf v}_n 0 \Longrightarrow \mathrm{\bf v}_n a : \mathrm{\bf v}_n C \longrightarrow \mathrm{\bf v}_n D$$
\[ \xymatrix@C=35pt@R=35pt{
\cdots \ar[r]^d & C^0 \ar[r]^{d^{n+1}} & C^{n+1} \ar[r]^d \ar[dl]_(.55){\hat{h}'} & C^{n+2} \ar[r]^{d^{n+1}} \ar[dl]|(.45){\hat{h}' d^n + d \hat{h}' d^{n-1} + \cdots + d^n \hat{h}'} & \cdots \\
\cdots \ar[r]^d & D^0 \ar[r]^{d^{n+1}} & D^{n+1} \ar[r]^d & D^{n+2} \ar[r]^{d^{n+1}} & \cdots ,
} \]
define the homotopy of $(n+3)$-complexes $\G_n(\hat{h}' , \hat{h}'')$ as:
\tiny
\[ \xymatrix@C=10pt@R=30pt{
\ar[r] & C^{n+2} \oplus X^2 \ar[r]^{\mtrx{d}{0}{0}{1}} \ar[dr]|{\mtrx{\hat{h}'}{0}{0}{0}} & C^{n+3} \oplus X^2 \ar[r]^(.65){\mtrx{d}{0}{0}{1}} \ar[dr]|{\mtrx{\hat{h}'}{0}{0}{0}} & \cdots \ar[r]^(.3){\mtrx{d}{0}{0}{1}} & C^{2n+3} \oplus X^2 \ar[r]^(.55){\mtrx{-d}{0}{x}{d^\star}} \ar[dr]|{\mtrx{-\hat{h}'}{0}{0}{\hat{h}_1}} & C^{2n+4} \oplus X^3 \ar[r]^{\mtrx{-1}{0}{x}{d}} \ar[dr]|{\mtrx{0}{0}{0}{\hat{h}_1}} & C^{2n+4} \oplus X^4 \ar[r]^(.67){\mtrx{d}{0}{0}{1}} & \\
\ar[r] & D^0 \oplus Y^0 \ar[r]_{\mtrx{d}{0}{0}{1}} & D^1 \oplus Y^0 \ar[r]_(.45){\mtrx{d}{0}{0}{1}} & D^2 \oplus Y^0 \ar[r]_(.55){\mtrx{d}{0}{0}{1}} & \cdots \ar[r]_(.55){\mtrx{-d}{0}{y}{d^\star}} & D^{n+2} \oplus Y^1 \ar[r]_{\mtrx{-1}{0}{y}{d}} & D^{n+2} \oplus Y^2 \ar[r]_(.67){\mtrx{d}{0}{0}{1}} & .
} \]
\normalsize
One must verify that $\G_n(\hat{h}' , \hat{h}'')$ is indeed a 2-morphism $\G_n(0, 0) \Rightarrow \G_n(a, f)$.
At degree $n + 1$, the summands $d^{n+2-j} \G_n(\hat{h}', \hat{h}'') d^j : \G_n(C, x)^{n+1} \to \G_n(D, y)^{n+1}$ are:
\begin{align*}
\mtrx{d}{0}{0}{1}^{n+1} \mtrx{-1}{0}{y}{d} \mtrx{-\hat{h}'}{0}{0}{\hat{h}_1} & = \mtrx{d^{n+1} \hat{h}'}{0}{-y \hat{h}'}{d \hat{h}_1} & (j = 0) \\
\mtrx{d}{0}{0}{1}^{n+1} \mtrx{0}{0}{0}{\hat{h}_1} \mtrx{-d}{0}{x}{d^\star} & = \mtrx{0}{0}{\hat{h}_1 x}{\hat{h}_1 d^\star} & (j = 1) \\
\mtrx{d}{0}{0}{1}^{n+2-j} \mtrx{\hat{h}'}{0}{0}{0} \mtrx{d}{0}{0}{1}^{j-2} \mtrx{-1}{0}{x}{d} \mtrx{-d}{0}{x}{d^\star} & = \mtrx{d^{n+2-j} \hat{h}' d^{j-1}}{0}{0}{0} & (2 \leq j \leq n+2)
\end{align*}
And their sum is
$$\mtrx{d^{n+1} \hat{h}'}{0}{-y \hat{h}'}{d \hat{h}_1} + \mtrx{0}{0}{\hat{h}_1 x}{\hat{h}_1 d^\star} + \sum_{i = 0}^n \mtrx{d^{n-i} \hat{h}' d^{i+1}}{0}{0}{0} = \mtrx{a}{0}{\hat{f}}{f_1} = \G_n(a, f)^{n + 1}$$
The verification is similar at other degrees.
\end{definition}

\begin{theorem}\label{FG}
The two functors 
\[ \xymatrix@C=35pt{
\S_{\ord{n}} \, {}_{\mathrm{\bf v}_n} \mkern-8mu \times_{\dom} \LArr(\S_{\ord{0}}) \ar@<.5ex>[r]^(.7){\G_n} \ar@{<-}@<-.5ex>[r]_(.7){\F_n} & \S_{\ord{n+1}}
} \]
are equimorphisms (see Definition \ref{equimorphism}) in the sestercategory (see Definition \ref{sestercategory}) of small sesquicategories.
More specifically, there are natural transformations
\[ \xymatrix@C=35pt{
\G_n \F_n \ar@<.5ex>[r]^(.5){\theta} \ar@{<-}@<-.5ex>[r]_(.5){\eta} & \Id_{\S_{\ord{n+1}}}
} \]
and lax natural transformations (see Definition \ref{LFun-ob})
\[ \xymatrix@C=35pt{
\F_n \G_n \ar@<.5ex>[r]^(.35){\epsilon} \ar@{<-}@<-.5ex>[r]_(.35){\zeta} & \Id_{\S_{\ord{n}} \, {}_{\mathrm{\bf v}_n} \mkern-8mu \times_{\dom} \LArr(\S_{\ord{0}})}
} \]
such that
$$\theta  \eta = \id_\Id, \: \eta \. \theta \simeq \id_{\G_n \F_n}$$
and
$$\epsilon \. \zeta = \id_\Id, \: \zeta \epsilon \simeq \id_{\F_n \G_n}$$
where $\simeq$ denotes homotopy equivalence, that is existence of invertible modifications.
\end{theorem}
\begin{proof}
The functor $\G_n \F_n$ maps an $(n+3)$-complex $C$ to
\tiny
\[ \xymatrix@C=12pt@R=20pt{
\cdots \ar[r] & C^0 \oplus C^0 \ar[r]^{\mtrx{d}{0}{0}{1}} & C^1 \oplus C^0 \ar[r]^(.6){\mtrx{d}{0}{0}{1}} & \cdots \ar[r]^(.3){\mtrx{d}{0}{0}{1}} & C^{n+1} \oplus C^0 \ar[r]^{\mtrx{-d^2}{0}{d}{d^{n+2}}} & C^{n+3} \oplus C^{n+2} \ar[r]^{\mtrx{-1}{0}{1}{d}} & C^{n+3} \oplus C^{n+3} \ar[r]^(.65){\mtrx{d}{0}{0}{1}} & \cdots
} \]
\normalsize
and morphisms are mapped degreewise.

Define the natural transformations 
\[ \xymatrix@C=35pt{
\G_n \F_n \ar@<.5ex>[r]^(.5){\theta} \ar@{<-}@<-.5ex>[r]_(.5){\eta} & \Id_{\S_{\ord{n+1}}}
} \]
on $C$ as follows
\tiny
\[ \xymatrix@C=12pt@R=25pt{
\cdots \ar[r] & C^0 \oplus C^0 \ar[r]^{\mtrx{d}{0}{0}{1}} & C^1 \oplus C^0 \ar[r]^(.6){\mtrx{d}{0}{0}{1}} & \cdots \ar[r]^(.3){\mtrx{d}{0}{0}{1}} & C^{n+1} \oplus C^0 \ar[r]^{\mtrx{-d^2}{0}{d}{d^{n+2}}} & C^{n+3} \oplus C^{n+2} \ar[r]^{\mtrx{-1}{0}{1}{d}} & C^{n+3} \oplus C^{n+3} \ar[r]^(.65){\mtrx{d}{0}{0}{1}} & \cdots \\
\cdots \ar[r] & C^0 \ar@{<-}@<.5ex>[u]^(.5){\rvec{1}{1}} \ar@<-.5ex>[u]_(.5){\cvec{1}{0}} \ar[r]_d & C^1 \ar@{<-}@<.5ex>[u]^(.5){\rvec{1}{d}} \ar@<-.5ex>[u]_(.5){\cvec{1}{0}} \ar[r]_d & \cdots \ar[r]_d & C^{n+1} \ar@{<-}@<.5ex>[u]^(.5){\rvec{1}{d^{n+1}}} \ar@<-.5ex>[u]_(.5){\cvec{1}{0}} \ar[r]_d & C^{n+2} \ar@{<-}@<.5ex>[u]^(.5){\rvec{0}{1}} \ar@<-.5ex>[u]_(.5){\cvec{-d}{1}} \ar[r]_d & C^{n+3} \ar@{<-}@<.5ex>[u]^(.5){\rvec{1}{1}} \ar@<-.5ex>[u]_(.5){\cvec{1}{0}} \ar[r]_d & \cdots
} \]
\normalsize
where $\theta_C$ is going downwards and $\eta_C$ is going upwards. They both commute with the differentials and they form natural transformations because morphisms are given degreewise. It is easy to see that
$$\theta \eta = \id_{\Id_{\S_{\ord{n+1}}}}.$$
The natural transformation $\eta \theta$ is not equal to, but is homotopy equivalent to $\id_{\G_n \F_n}$. To see this, we compute the morphism $\id_{\G_n \F_n C} - \eta_C \theta_C$
\tiny
\[ \xymatrix@C=12pt@R=25pt{
\cdots \ar[r] & C^0 \oplus C^0 \ar[d]|(.45){\mtrx{0}{-1}{0}{1}} \ar[r]^{\mtrx{d}{0}{0}{1}} & C^1 \oplus C^0 \ar[d]|(.45){\mtrx{0}{-d}{0}{1}} \ar[r]^(.6){\mtrx{d}{0}{0}{1}} & \cdots \ar[r]^(.3){\mtrx{d}{0}{0}{1}} & C^{n+1} \oplus C^0 \ar[d]|(.45){\mtrx{0}{-d^{n+1}}{0}{1}} \ar[r]^{\mtrx{d^2}{0}{d}{-d^{n+2}}} & C^{n+3} \oplus C^{n+2} \ar[d]|(.45){\mtrx{1}{d}{0}{0}} \ar[r]^{\mtrx{1}{0}{1}{-d}} & C^{n+3} \oplus C^{n+3} \ar[d]|(.45){\mtrx{0}{-1}{0}{1}} \ar[r]^(.65){\mtrx{d}{0}{0}{1}} & \cdots \\
\cdots \ar[r] & C^0 \oplus C^0 \ar[r]_{\mtrx{d}{0}{0}{1}} & C^1 \oplus C^0 \ar[r]_(.6){\mtrx{d}{0}{0}{1}} & \cdots \ar[r]_(.3){\mtrx{d}{0}{0}{1}} & C^{n+1} \oplus C^0 \ar[r]_{\mtrx{d^2}{0}{d}{-d^{n+2}}} & C^{n+3} \oplus C^{n+2} \ar[r]_{\mtrx{1}{0}{1}{-d}} & C^{n+3} \oplus C^{n+3} \ar[r]_(.65){\mtrx{d}{0}{0}{1}} & \cdots
} \]
\normalsize
The homotopy $h_C$ is given by
\[
\hat{h}_C^i = \begin{cases} \mtrx{0}{1}{0}{0} &\mbox{if } i \equiv -1 \pmod{n+3}\\
0 & \mbox{otherwise.} \end{cases}
\]
To verify that
$$\id_{\G_n \F_n C} - \eta_C \theta_C = \sum_{j = 1}^{n+3} d^{n+3-j} \. \hat{h}_C \. d^{j-1}$$
at a degree congruent to $i$ modulo $n+3$, where $0 \leq i \leq n+1$, one checks that
$$\mtrx{d}{0}{0}{1}^i \mtrx{-1}{0}{1}{d} \mtrx{0}{1}{0}{0} \mtrx{d}{0}{0}{1}^{n+1-i} = \mtrx{0}{-d^i}{0}{1}$$
and at a degree congruent to $n+2$ modulo $n+3$, one checks that
$$\mtrx{0}{1}{0}{0} \mtrx{d}{0}{0}{1}^{n+1} \mtrx{-1}{0}{1}{d} = \mtrx{1}{d}{0}{0}.$$

The functor $\F_n \G_n$ maps an object $(C, x : \mathrm{\bf v}_n C \to X)$ of $\S_{\ord{n}} \, {}_{\mathrm{\bf v}_n} \mkern-8mu \times_{\dom} \LArr(\S_{\ord{0}})$ to the pair consisting of the $(n+2)$-complex $\mathrm{\bf d}_{n+1} \G_n (C, x)$, namely,
\[ \xymatrix@C=15pt@R=35pt{
\cdots \ar[r] & C^0 \oplus X^0 \ar[r]^{\mtrx{d}{0}{0}{1}} & C^1 \oplus X^0 \ar[r]^(.6){\mtrx{d}{0}{0}{1}} & \cdots \ar[r]^(.3){\mtrx{d}{0}{0}{1}} & C^{n+1} \oplus X^0 \ar[r]^{\mtrx{d}{0}{0}{0}} & C^{n+2} \oplus X^2 \ar[r]^(.65){\mtrx{d}{0}{0}{1}} & \cdots
} \]
and the morphism of $2$-complexes $\tau_{\G_n (C, x)} : \mathrm{\bf v}_{n} \G_n (C, x) \to \mathrm{\bf v}_{n+1} \G_n (C, x)$ given as follows:
\[ \xymatrix@C=35pt@R=35pt{
\cdots \ar[r] & C^0 \oplus X^0 \ar[r]^(.45){\mtrx{d^{n+1}}{0}{0}{1}} \ar[d]|{\mtrx{1}{0}{0}{1}} & C^{n+1} \oplus X^0 \ar[r]^{\mtrx{d}{0}{0}{0}} \ar[d]|{\mtrx{-d}{0}{x}{d^\star}} & C^{n+2} \oplus X^2 \ar[r]^(.65){\mtrx{d^{n+1}}{0}{0}{1}} \ar[d]|{\mtrx{1}{0}{0}{1}} & \cdots \\
\cdots \ar[r] & C^0 \oplus X^0 \ar[r]_(.45){\mtrx{0}{0}{xd^{n+1}}{d^\star}} & C^{n+2} \oplus X^1 \ar[r]_{\mtrx{-1}{0}{x}{d}} & C^{n+2} \oplus X^2 \ar[r]_(.65){\mtrx{0}{0}{xd^{n+1}}{d^\star}} & \cdots
} \]
and $\F_n \G_n$ maps a morphism $(a, f) : (C, x) \to (D, y)$ to the pair consisting of the morphism of $(n+2)$-complexes $\mathrm{\bf d}_{n+1} \G_n (a, f)$
\[ \xymatrix@C=14pt@R=35pt{
\cdots \ar[r] & C^0 \oplus X^0 \ar[r]^{\mtrx{d}{0}{0}{1}} \ar[d]|{\mtrx{a}{0}{\hat{f}d^{n+1}}{f_1}} & C^1 \oplus X^0 \ar[r]^(.6){\mtrx{d}{0}{0}{1}} \ar[d]|{\mtrx{a}{0}{\hat{f}d^n}{f_1}} & \cdots \ar[r]^(.3){\mtrx{d}{0}{0}{1}} & C^{n+1} \oplus X^0 \ar[r]^{\mtrx{d}{0}{0}{0}} \ar[d]|{\mtrx{a}{0}{\hat{f}}{f_1}} & C^{n+2} \oplus X^2 \ar[r]^(.65){\mtrx{d}{0}{0}{1}} \ar[d]|{\mtrx{a}{0}{\hat{f}d^{n+1}}{f_1}} & \cdots \\
\cdots \ar[r] & D^0 \oplus Y^0 \ar[r]_{\mtrx{d}{0}{0}{1}} & D^1 \oplus Y^0 \ar[r]_(.6){\mtrx{d}{0}{0}{1}} & \cdots \ar[r]_(.3){\mtrx{d}{0}{0}{1}} & D^{n+1} \oplus Y^0 \ar[r]_{\mtrx{d}{0}{0}{0}} & D^{n+2} \oplus Y^2 \ar[r]_(.65){\mtrx{d}{0}{0}{1}} & \cdots
} \]
and the commuting diagram of $2$-complexes $\tau_{\G_n (a, f)} {\kern -1.5pt} = {\kern -1.5pt} (\mathrm{\bf v}_{n} \G_n (a, f), \mathrm{\bf v}_{n+1} \G_n (a, f), 0)$
\[ \xymatrix@C=35pt@R=35pt{
\mathrm{\bf v}_{n} \G_n (C, x) \ar[d]_{\tau_{\G_n (C, x)}} \ar[r]^{\mathrm{\bf v}_{n} \G_n (a, f)} 
& \mathrm{\bf v}_{n} \G_n (D, y) \ar[d]^{\tau_{\G_n (D, y)}} \\
\mathrm{\bf v}_{n+1} \G_n (C, x) \ar[r]_{\mathrm{\bf v}_{n+1} \G_n (a, f)} &
\mathrm{\bf v}_{n+1} \G_n (D, y) 
} \]
where $\mathrm{\bf v}_{n} \G_n (a, f)$ is given by
\[ \xymatrix@C=35pt@R=35pt{
\cdots \ar[r] & C^0 \oplus X^0 \ar[r]^{\mtrx{d^{n+1}}{0}{0}{1}} \ar[d]|{\mtrx{a}{0}{\hat{f}d^{n+1}}{f_1}} & C^{n+1} \oplus X^0 \ar[r]^{\mtrx{d}{0}{0}{0}} \ar[d]|{\mtrx{a}{0}{\hat{f}}{f_1}} & C^{n+2} \oplus X^2 \ar[r]^(.6){\mtrx{d^{n+1}}{0}{0}{1}} \ar[d]|{\mtrx{a}{0}{\hat{f}d^{n+1}}{f_1}} & \cdots \\
\cdots \ar[r] & D^0 \oplus Y^0 \ar[r]_{\mtrx{d^{n+1}}{0}{0}{1}} & D^{n+1} \oplus Y^0 \ar[r]_{\mtrx{d}{0}{0}{0}} & D^{n+2} \oplus Y^2 \ar[r]_(.6){\mtrx{d^{n+1}}{0}{0}{1}} & \cdots
} \]
and $\mathrm{\bf v}_{n+1} \G_n (a, f)$ is given by
\[ \xymatrix@C=35pt@R=35pt{
\cdots \ar[r] & C^0 \oplus X^0 \ar[r]^{\mtrx{0}{0}{xd^{n+1}}{d^\star}} \ar[d]|{\mtrx{a}{0}{\hat{f}d^{n+1}}{f_1}} & C^{n+2} \oplus X^1 \ar[r]^{\mtrx{-1}{0}{x}{d}} \ar[d]|{\mtrx{a}{0}{\hat{f}}{f_1}} & C^{n+2} \oplus X^2 \ar[r]^(.6){\mtrx{0}{0}{xd^{n+1}}{d^\star}} \ar[d]|{\mtrx{a}{0}{\hat{f}d^{n+1}}{f_1}} & \cdots \\
\cdots \ar[r] & D^0 \oplus Y^0 \ar[r]_{\mtrx{0}{0}{yd^{n+1}}{d^\star}} & D^{n+2} \oplus Y^1 \ar[r]_{\mtrx{-1}{0}{y}{d}} & D^{n+2} \oplus Y^2 \ar[r]_(.6){\mtrx{0}{0}{yd^{n+1}}{d^\star}} & \cdots .
} \]

Now, it remains to define two lax natural transformations
\[ \xymatrix@C=35pt{
\F_n \G_n \ar@<.5ex>[r]^(.35){\epsilon} \ar@{<-}@<-.5ex>[r]_(.35){\zeta} & \Id_{\S_{\ord{n}} \, {}_{\mathrm{\bf v}_n} \mkern-8mu \times_{\dom} \LArr(\S_{\ord{0}})}
} \]
such that $\epsilon \. \zeta = \id_{\Id_{\S_{\ord{n}} \, {}_{\mathrm{\bf v}_n} \mkern-8mu \times_{\dom} \LArr(\S_{\ord{0}})}}$ and $\zeta \epsilon$ is homotopy equivalent to $\id_{\F_n \G_n}$.
Each of the lax natural transformations $\epsilon$ and $\zeta$ is a pair of lax natural transformations
$$\epsilon = (\epsilon', \epsilon''), \: \zeta = (\zeta', \zeta'')$$
so that the components match on $\mathrm{\bf v}_n$ and $\dom$.
Recall that these lax natural transformations are, by Definition \ref{LFun-ob}, functors
\begin{align*}
\epsilon', \zeta' : \S_{\ord{n}} \, {}_{\mathrm{\bf v}_n} \mkern-8mu \times_{\dom} \LArr(\S_{\ord{0}}) & \longrightarrow \LArr(\S_{\ord{n}}), \\
\epsilon'', \zeta'' : \S_{\ord{n}} \, {}_{\mathrm{\bf v}_n} \mkern-8mu \times_{\dom} \LArr(\S_{\ord{0}}) & \longrightarrow \LArr(\LArr(\S_{\ord{0}})).
\end{align*}

Define $\epsilon'$ and $\zeta'$ on objects $(C, x)$ as follows:
\[\begin{array}{rr}
\epsilon'_{(C, x)} := & \rvec{1}{0}  : \mathrm{\bf d}_{n+1} \G_n (C, x) \longrightarrow C, \\
\zeta'_{(C, x)} := & \cvec{1}{0}  : C \longrightarrow \mathrm{\bf d}_{n+1} \G_n (C, x),
\end{array}\]
that is
\[ \xymatrix@C=15pt@R=35pt{
\cdots \ar[r] & C^0 \oplus X^0 \ar[r]^{\mtrx{d}{0}{0}{1}} & C^1 \oplus X^0 \ar[r]^(.6){\mtrx{d}{0}{0}{1}} & \cdots \ar[r]^(.3){\mtrx{d}{0}{0}{1}} & C^{n+1} \oplus X^0 \ar[r]^{\mtrx{d}{0}{0}{0}} & C^{n+2} \oplus X^2 \ar[r]^(.65){\mtrx{d}{0}{0}{1}} & \cdots \\
\cdots \ar[r] & C^0 \ar@{<-}@<.5ex>[u]^(.5){\rvec{1}{0}} \ar@<-.5ex>[u]_(.5){\cvec{1}{0}} \ar[r]_d & C^1 \ar@{<-}@<.5ex>[u]^(.5){\rvec{1}{0}} \ar@<-.5ex>[u]_(.5){\cvec{1}{0}} \ar[r]_d & \cdots \ar[r]_d & C^{n+1} \ar@{<-}@<.5ex>[u]^(.5){\rvec{1}{0}} \ar@<-.5ex>[u]_(.5){\cvec{1}{0}} \ar[r]_d & C^{n+2} \ar@{<-}@<.5ex>[u]^(.5){\rvec{1}{0}} \ar@<-.5ex>[u]_(.5){\cvec{1}{0}} \ar[r]_d & \cdots .
} \]
On a morphisms $(a, f) : (C, x) \to (D, y)$, $\epsilon'$ and $\zeta'$ are defined as follows:
\begin{align*}
\epsilon'_{(a, f)} &= (\mathrm{\bf d}_{n+1} \G_n (a, f), a, \hat{\epsilon}'_{(a, f)}), \\
\zeta'_{(a, f)} &= (a, \mathrm{\bf d}_{n+1} \G_n (a, f), \hat{\zeta}'_{(a, f)}),
\end{align*}
that is
\[ \xymatrix@C=40pt@R=45pt{
C \ar[d]_a \ar[r]^(.35){\zeta'_{(C, x)}} \ar@{.>}[dr]_{\hat{\zeta}'_{(a, f)}} \xtwocell[dr]{}\omit{^} & \mathrm{\bf d}_{n+1} \G_n (C, x) \ar[d]|(.42){\mathrm{\bf d}_{n+1} \G_n (a, f)} \ar[r]^(.65){\epsilon'_{(C, x)}} \ar@{.>}[dr]_{\hat{\epsilon}'_{(a, f)}} \xtwocell[dr]{}\omit{^} & C \ar[d]^a \\
D \ar[r]_(.35){\zeta'_{(D, y)}} & \mathrm{\bf d}_{n+1} \G_n (D, y) \ar[r]_(.65){\epsilon'_{(D, y)}} & D 
} \]
where $\hat{\epsilon}'_{(a, f)}$ and $\hat{\zeta}'_{(a, f)}$ are defined as follows:
\begin{align*}
\hat{\epsilon}'_{(a, f)} &:= 0, \\
\left(\hat{\zeta}'_{(a, f)}\right)^i &:= \begin{cases} \cvec{0}{\hat{f}} &\mbox{if } i \equiv 0 \pmod{n+2} ;\\
\cvec{0}{\vphantom{\hat{f}}0} & \mbox{otherwise.} \end{cases}
\end{align*}
The lax natural transformation $\hat{\epsilon}'$ is hence a natural transformation.

The lax natural transformations $\epsilon''$ and $\zeta''$ map an object $(C, x)$ to
\[\xymatrix@C=35pt@R=35pt{
\mathrm{\bf v}_{n} \G_n (C, x) \ar[d]_{\tau_{\G_n (C, x)}} \ar[r]^(.6){\epsilon''_{(C, x) 0}} \ar@{.>}[dr]^(.6){\hat{\epsilon}''_{(C, x)}} \xtwocell[dr]{}\omit{_} & \mathrm{\bf v}_n C \ar[d]^x \\
\mathrm{\bf v}_{n+1} \G_n (C, x) \ar[r]_(.6){\epsilon''_{(C, x) 1}} & X }
\xymatrix@C=0pt@R=15pt{ \\ \mbox{ and }}
\xymatrix@C=35pt@R=35pt{
\mathrm{\bf v}_n C \ar[d]_x \ar[r]^(.4){\zeta''_{(C, x) 0}} \ar@{.>}[dr]^(.6){\hat{\zeta}''_{(C, x)}} \xtwocell[dr]{}\omit{_} & \mathrm{\bf v}_{n} \G_n (C, x) \ar[d]^{\tau_{\G_n (C, x)}} \\
X \ar[r]_(.4){\zeta''_{(C, x) 1}} & \mathrm{\bf v}_{n+1} \G_n (C, x) \\
}\]
where
$$\epsilon''_{(C, x) 0}  := \mathrm{\bf v}_{n} \epsilon'_{(C, x)} = \rvec{1}{0} \mbox{ and } \zeta''_{(C, x) 0} := \mathrm{\bf v}_{n} \zeta'_{(C, x)} = \cvec{1}{0}$$
by necessity, $\epsilon''_{(C, x) 1}$ is given by
\[ \xymatrix@C=30pt@R=35pt{
\cdots \ar[r] & C^0 \oplus X^0 \ar[r]^(.47){\mtrx{0}{0}{xd^{n+1}}{d^\star}} & C^{n+2} \oplus X^1 \ar[r]^{\mtrx{-1}{0}{x}{d}} & C^{n+2} \oplus X^2 \ar[r]^(.65){\mtrx{0}{0}{xd^{n+1}}{d^\star}} & \cdots \\
\cdots \ar[r] & X^0 \ar[r]^{d^\star} \ar@{<-}[u]_{\rvec{x}{1}} & X^1 \ar[r]^d \ar@{<-}[u]_{\rvec{0}{1}} & X^2 \ar[r]^{d^\star} \ar@{<-}[u]_{\rvec{x}{1}} & \cdots ,
} \]
$\zeta''_{(C, x) 1}$ is given by
\[ \xymatrix@C=30pt@R=35pt{
\cdots \ar[r] & X^0 \ar[r]^{d^\star} \ar[d]^{\cvec{0}{1}} & X^1 \ar[r]^d \ar[d]^{\cvec{0}{1}} & X^2 \ar[r]^{d^\star} \ar[d]^{\cvec{0}{1}} & \cdots \\
\cdots \ar[r] & C^0 \oplus X^0 \ar[r]_(.45){\mtrx{0}{0}{xd^{n+1}}{{d^\star}}} & C^{n+2} \oplus X^1 \ar[r]_{\mtrx{-1}{0}{x}{d}} & C^{n+2} \oplus X^2 \ar[r]_(.65){\mtrx{0}{0}{xd^{n+1}}{{d^\star}}} & \cdots ,
} \]
$\hat{\epsilon}''_{(C, x)}$ is given by
\[ \xymatrix@C=30pt@R=35pt{
\cdots \ar[r] & C^0 \oplus X^0 \ar[r]^(.47){\mtrx{d^{n+1}}{0}{0}{1}} & C^{n+1} \oplus X^0 \ar[r]^{\mtrx{d}{0}{0}{0}} \ar[dl]_(.55){\rvec{0}{1}} & C^{n+2} \oplus X^2 \ar[r]^(.6){\mtrx{d^{n+1}}{0}{0}{1}} \ar[dl]_(.55){\rvec{0}{0}} & \cdots \\
\cdots \ar[r] & X^0 \ar[r]^{d^\star} & X^1 \ar[r]^d & X^2 \ar[r]^{d^\star} & \cdots 
} \]
and $\hat{\zeta}''_{(C, x)}$ is given by
\[ \xymatrix@C=30pt@R=35pt{
\cdots \ar[r] & C^0 \ar[r]^{d^{n+1}} & C^{n+1} \ar[r]^d \ar[dl]_(.6){\cvec{0}{0}} & C^{n+2} \ar[r]^{d^{n+1}} \ar[dl]_(.6){\cvec{1}{0}} & \cdots \\
\cdots \ar[r] & C^0 \oplus X^0 \ar[r]_(.45){\mtrx{0}{0}{xd^{n+1}}{{d^\star}}} & C^{n+2} \oplus X^1 \ar[r]_{\mtrx{-1}{0}{x}{d}} & C^{n+2} \oplus X^2 \ar[r]_(.65){\mtrx{0}{0}{xd^{n+1}}{{d^\star}}} & \cdots .
} \]

On a morphism $(a, f) : (C, x) \to (D, y)$, $\epsilon''_{(a, f)}$ is a diagram
\[ \xymatrix@C=40pt@R=40pt{
\tau_{\G_n (C, x)} \ar[r]^{\tau_{\G_n (a, f)}} \ar[d]_{\epsilon''_{(C, x)}} \ar@{.>}[dr]^(.45){\hat{\epsilon}''_{(a, f)}} \xtwocell[dr]{}\omit{_} & \tau_{\G_n (D, y)} \ar[d]^{\epsilon''_{(D, y)}} \\
x \ar[r]_{f} & y
} \]
or, more explicitly,
\[ \xymatrix@C=35pt@R=35pt{
\mathrm{\bf v}_{n} \G_n (C, x) \ar[dddd]_{\tau_{\G_n (C, x)}} \ar[rrrr]^{\mathrm{\bf v}_{n} \G_n (a, f)} \ar[dr]^(.6){\epsilon''_{(C, x) 0}} \ar@{.>}[drrr]^(.6){\hat{\epsilon}''_{(a, f) 0}} \xtwocell[drrr]{}\omit{_} \ar@{.>}[dddr]_(.6){\hat{\epsilon}''_{(C, x)}} \xtwocell[dddr]{}\omit{_} & & & & \mathrm{\bf v}_{n} \G_n (D, y) \ar[dl]_{\epsilon''_{(D, y) 0}} \ar[dddd]^{\tau_{\G_n (D, y)}} \ar@{.>}[dddl]^(.6){\hat{\epsilon}''_{(D, y)}} \xtwocell[dddl]{}\omit{^} \\
& \mathrm{\bf v}_{n} C \ar[dd]^{x} \ar[rr]_{\mathrm{\bf v}_{n} a} \ar@{.>}[ddrr]^(.55){\hat{f}} \xtwocell[ddrr]{}\omit{_} & & \mathrm{\bf v}_{n} D \ar[dd]_{y} &\\
&&&&\\
& X \ar[rr]^{f_1} & & Y &\\
\mathrm{\bf v}_{n+1} \G_n (C, x) \ar[ur]^(.6){\epsilon''_{(C, x) 1}} \ar[rrrr]_{\mathrm{\bf v}_{n+1} \G_n (a, f)} \ar@{.>}[urrr]_(.6){\hat{\epsilon}''_{(a, f) 1}} \xtwocell[urrr]{}\omit{^} & & & & \mathrm{\bf v}_{n+1} \G_n (D, y) \ar[ul]_(.6){\epsilon''_{(D, y) 1}} \. ,\\
}\]
where
$$\hat{\epsilon}''_{(a, f) 0} := 0$$
and $\hat{\epsilon}''_{(a, f) 1}$ is given by
\[ \xymatrix@C=30pt@R=35pt{
\cdots \ar[r] & C^0 \oplus X^0 \ar[r]^(.47){\mtrx{0}{0}{xd^{n+1}}{{d^\star}}} & C^{n+2} \oplus X^1 \ar[r]^{\mtrx{-1}{0}{x}{d}} \ar[dl]_(.55){\rvec{0}{0}}  & C^{n+2} \oplus X^2 \ar[r]^(.65){\mtrx{0}{0}{xd^{n+1}}{{d^\star}}} \ar[dl]_(.55){\rvec{\hat{f}}{0}} & \cdots \\
\cdots \ar[r] & Y^0 \ar[r]^{d^\star} & Y^1 \ar[r]^d & Y^2 \ar[r]^{d^\star} & \cdots .
} \]

Likewise, $\zeta''_{(a, f)}$ is a diagram
\[ \xymatrix@C=40pt@R=40pt{
\tau_{\G_n (C, x)} \ar[r]^{\tau_{\G_n (a, f)}} \ar@{<-}[d]_{\zeta''_{(C, x)}} & \tau_{\G_n (D, y)} \ar@{<-}[d]^{\zeta''_{(D, y)}} \\
x \ar[r]_{f} \ar@{.>}[ur]_(.6){\hat{\zeta}''_{(a, f)}} \xtwocell[ur]{}\omit{^} & y
} \]
or more explicitly
\[ \xymatrix@C=35pt@R=35pt{
\mathrm{\bf v}_{n} \G_n (C, x) \ar[dddd]_{\tau_{\G_n (C, x)}} \ar[rrrr]^{\mathrm{\bf v}_{n} \G_n (a, f)} \ar@{<-}[dr]_(.6){\zeta''_{(C, x) 0}}
& & & & \mathrm{\bf v}_{n} \G_n (D, y) \ar@{<-}[dl]^(.6){\zeta''_{(D, y) 0}} \ar[dddd]^{\tau_{\G_n (D, y)}}  \\
& \mathrm{\bf v}_{n} C \ar[dd]^{x} \ar[rr]_{\mathrm{\bf v}_{n} a}
\ar@{.>}[urrr]^(.4){\hat{\zeta}''_{(a, f) 0}} \xtwocell[urrr]{}\omit{^}
\ar@{.>}[dddl]_(.4){\hat{\zeta}''_{(C, x)}} \xtwocell[dddl]{}\omit{^}
\ar@{.>}[ddrr]^(.55){\hat{f}} \xtwocell[ddrr]{}\omit{_}
& & \mathrm{\bf v}_{n} D \ar[dd]_{y}
\ar@{.>}[dddr]^(.4){\hat{\zeta}''_{(D, y)}} \xtwocell[dddr]{}\omit{_} &\\
&&&&\\
& X \ar[rr]^{f_1} \ar@{.>}[rrrd]_(.4){\hat{\zeta}''_{(a, f) 1}} \xtwocell[rrrd]{}\omit{_} & & Y &\\
\mathrm{\bf v}_{n+1} \G_n (C, x) \ar@{<-}[ur]^(.65){\zeta''_{(C, x) 1}} \ar[rrrr]_{\mathrm{\bf v}_{n+1} \G_n (a, f)} & & & & \mathrm{\bf v}_{n+1} \G_n (D, y) \ar@{<-}[ul]_(.65){\zeta''_{(D, y) 1}} \. ,\\
}\]
where $\hat{\zeta}''_{(a, f)0}$ is given by
\[ \xymatrix@C=35pt@R=35pt{
\cdots \ar[r] & C^0 \ar[r]^{d^{n+1}} & C^{n+1} \ar[r]^{d} \ar[dl]_(.6){\cvec{0}{\hat{f}}} & C^{n+2} \ar[r]^{d^{n+1}} \ar[dl]_(.6){\cvec{0}{0\vphantom{\hat{f}}}} & \cdots \\
\cdots \ar[r] & D^0 \oplus Y^0 \ar[r]_{\mtrx{d^{n+1}}{0}{0}{1}} & D^{n+1} \oplus Y^0 \ar[r]_{\mtrx{d}{0}{0}{0}} & D^{n+2} \oplus Y^2 \ar[r]_(.6){\mtrx{d^{n+1}}{0}{0}{1}} & \cdots ,
} \]
and 
$$\hat{\zeta}''_{(a, f)1} := 0.$$

Now, one has to verify that
\begin{itemize}
\item $\epsilon', \zeta'$ are functors $\S_{\ord{n}} \, {}_{\mathrm{\bf v}_n} \mkern-8mu \times_{\dom} \LArr(\S_{\ord{0}}) \to \LArr(\S_{\ord{n}})$,
\item $\epsilon'', \zeta''$ are functors $\S_{\ord{n}} \, {}_{\mathrm{\bf v}_n} \mkern-8mu \times_{\dom} \LArr(\S_{\ord{0}}) \to \LArr(\LArr(\S_{\ord{0}}))$,
\item $\LArr(\mathrm{\bf v}_n) \epsilon' = \LArr(\dom) \epsilon''$ and $\LArr(\mathrm{\bf v}_n) \zeta' = \LArr(\dom) \zeta''$,
\item $\epsilon' \zeta' = \id_{\Id_{\S_{\ord{n}}}}$ and $\epsilon'' \zeta'' = \id_{\Id_{\LArr(\S_{\ord{0}})}}$,
\item $\zeta' \epsilon' \simeq \id_{\mathrm{\bf v}_{n} \G_n}$ and $\zeta'' \epsilon'' \simeq \id_{\tau \G_n}$.
\end{itemize}
\end{proof}

\begin{proposition}\label{FG_adj}
The tuple
$$\left( \F_n, \G_n, \epsilon, \eta \right)$$
constitute an adjoint equivalence (see Definition \ref{adj_eq}).
\end{proposition}
\begin{proof}
We have already shown equivalence, it remains to prove the triangular identities
$$\epsilon_{\F_n} \! (\F_n \eta) = \id_{\F_n} , \: (\G_n \epsilon) \eta_{G_n} \!\! = \id_{\G_n}.$$

To show the first equality, let $C$ be an object of $\S_{\ord{n+1}}$, that is an $n+3$ complex.
The left hand side is the composition
$$\xymatrix@C=40pt{\F_n C \ar[r]^(.4){\F_n \eta_C} & \F_n \G_n \F_n C \ar[r]^(.55){\epsilon_{\F_n C}} & \F_n C}.$$
Its first component is a morphism between 2-complexes
$$\xymatrix@C=40pt{\mathrm{\bf d}_{n+1} C \ar[r]^(.4){\mathrm{\bf d}_{n+1} \eta_C} & \mathrm{\bf d}_{n+1} \G_n \F_n C \ar[r]^(.55){\epsilon_{\F_n C}'} & \mathrm{\bf d}_{n+1} C},$$
which is equal to
$$\rvec{1}{0} \cvec{1}{0} = 1.$$
On morphisms, $\eta$ is a natural transformation and $\hat{\epsilon}' = 0$.
Its second component is in $\LArr(\S_{\ord{n}})$
\[ \xymatrix@C=40pt@R=45pt{
\mathrm{\bf v}_{n} C \ar[d]_{\tau_C} \ar[r]^(.4){\mathrm{\bf v}_{n} \eta_C} & \mathrm{\bf v}_{n} \G_n \F_n C \ar[d]_{\tau_{\G_n \F_n C}} \ar[r]^(.6){\epsilon''_{(\F_n C) 0}} \ar@{.>}[dr]^(.6){\hat{\epsilon}''_{\F_n C}} \xtwocell[dr]{}\omit{_} & \mathrm{\bf v}_{n} C \ar[d]^{\tau_C} \\
\mathrm{\bf v}_{n+1} C \ar[r]_(.4){\mathrm{\bf v}_{n+1} \eta_C} & \mathrm{\bf v}_{n+1} \G_n \F_n C \ar[r]_(.6){\epsilon''_{(\F_n C) 1}} & \mathrm{\bf v}_{n+1} C .
} \]
We have already checked the upper row because the $n$:th vertex of the first component equals the domain of the second. The lower row is as follows
\[ \xymatrix@C=30pt@R=35pt{
\cdots \ar[r] & C^0 \ar[r]^{d^{n+2}} \ar[d]_{\cvec{1}{0}} & C^{n+2} \ar[r]^d \ar[d]^{\cvec{-d}{1}} & C^{n+3} \ar[r] \ar[d]^{\cvec{1}{0}} & \cdots \\
\cdots \ar[r] & C^0 \oplus C^0 \ar[r]^(.47){\mtrx{0}{0}{d^{n+2}}{d^{n+2}}} & C^{n+3} \oplus C^{n+2} \ar[r]^{\mtrx{-1}{0}{1}{d}} & C^{n+3} \oplus C^{n+3} \ar[r] & \cdots \\
\cdots \ar[r] & C^0 \ar[r]^{d^{n+2}} \ar@{<-}[u]^{\rvec{1}{1}} & C^{n+2} \ar[r]^d \ar@{<-}[u]_{\rvec{0}{1}} & C^{n+3} \ar[r] \ar@{<-}[u]_{\rvec{1}{1}} & \cdots 
} \]
which is identity on $\mathrm{\bf v}_{n+1} C$. And for the 2-morphism 
$$\hat{\epsilon}''_{\F_n C} \mathrm{\bf v}_{n} \eta_C = \rvec{0}{*} \cvec{1}{0} = 0 .$$
On morphisms $f : C \to D$, $\eta$ is a natural transformation, $\hat{\epsilon}''_{(\F_n f)0} = 0$ by definition and $\hat{\epsilon}''_{(\F_n f)1} = 0$ because the diagram
\[ \xymatrix@C=35pt@R=35pt{
\mathrm{\bf v}_{n} C \ar[d]_{\tau_{C}} \ar[r]^{\mathrm{\bf v}_{n} f} 
& \mathrm{\bf v}_{n} D \ar[d]^{\tau_{D}} \\
\mathrm{\bf v}_{n+1} C \ar[r]_{\mathrm{\bf v}_{n+1} f} &
\mathrm{\bf v}_{n+1} D
} \]
commutes. This proves the first equality, the second one can be checked in a similar way.
\end{proof}

\begin{theorem}\label{th_filler}
The functor $\Spine_{\ord{n}}$ has a weak inverse.
\end{theorem}
The theorem is proven by explicit construction of the weak inverse, which we denote $\Filler_{\ord{n}}$
\[ \xymatrix{
\S_{\ord{n}} \ar@<.5ex>[rr]^(.43){\Spine_{\ord{n}}} & &
\ar@<.5ex>[ll]^(.57){\Filler_{\ord{n}}} \LN(\S_{\ord{0}})_{\ord{n}} .
} \]
\begin{proof}
Fix $n$. Regard the category $\LN(\S_{\ord{0}})_{\ord{n}}$ as the limit of
\small
\[ \xymatrix@C=10pt@R=10pt{
\LArr(\S_{\ord{0}})_{1} \ar[dr]_{\cod} & &  \ar[dl]^{\dom} \LArr(\S_{\ord{0}})_{2} \ar@{}[dr]|(.52){\mbox{.}} \ar@{}[dr]|(.6){\mbox{.}} \ar@{}[dr]|(.68){\mbox{.}} & & \ar[dl]^{\dom} \LArr(\S_{\ord{0}})_{n} \\
& \S_{\ord{0}} & & \S_{\ord{0}} & & .
} \]
\normalsize
Introducing two $\S_{\ord{0}}$ at the left in the following way
\small
\[ \xymatrix@C=10pt@R=10pt{
\S_{\ord{0}} \ar[dr]_{\mathrm{\bf v}_0 = \Id \:\:\:} & & \ar[dl]^{\dom} \LArr(\S_{\ord{0}})_{1} \ar[dr]_{\cod} & &  \ar[dl]^{\dom} \LArr(\S_{\ord{0}})_{2} \ar@{}[dr]|(.52){\mbox{.}} \ar@{}[dr]|(.6){\mbox{.}} \ar@{}[dr]|(.68){\mbox{.}} & & \ar[dl]^{\dom} \LArr(\S_{\ord{0}})_{n} \\
& \S_{\ord{0}} & & \S_{\ord{0}} & & \S_{\ord{0}} &
} \]
\normalsize
does not change its limit. To see this, flip the direction of the identity morphism and observe that the subdiagram without the $\S_{\ord{0}}$ objects is initial (dual of final, see \cite[IX.3]{Mac}).

The idea of the proof is to utilize Theorem \ref{FG} to construct an equimorphism between the limit of
\small
\[ \xymatrix@C=10pt@R=10pt{
\S_{\ord{i}} \ar[dr]_{\mathrm{\bf v}_{i}} & & \ar[dl]^{\dom} \LArr(\S_{\ord{0}})_{i+1} \ar[dr]_{\cod} & &  \ar[dl]^{\dom} \LArr(\S_{\ord{0}})_{i+2} \ar@{}[dr]|(.52){\mbox{.}} \ar@{}[dr]|(.6){\mbox{.}} \ar@{}[dr]|(.68){\mbox{.}} & & \ar[dl]^{\dom} \LArr(\S_{\ord{0}})_{n} \\
& \S_{\ord{0}} & & \S_{\ord{0}} & & \S_{\ord{0}} 
} \]
\normalsize
denoted $\li_i$, and the limit of
\small
\[ \xymatrix@C=10pt@R=10pt{
\S_{\ord{i+1}} \ar[dr]_{\mathrm{\bf v}_{i+1}} & & \ar[dl]^{\dom} \LArr(\S_{\ord{0}})_{i+2} \ar[dr]_{\cod} & &  \ar[dl]^{\dom} \LArr(\S_{\ord{0}})_{i+3} \ar@{}[dr]|(.52){\mbox{.}} \ar@{}[dr]|(.6){\mbox{.}} \ar@{}[dr]|(.68){\mbox{.}} & & \ar[dl]^{\dom} \LArr(\S_{\ord{0}})_{n} \\
& \S_{\ord{0}} & & \S_{\ord{0}} & & \S_{\ord{0}} & &
} \]
\normalsize
denoted $\li_{i+1}$ and then to compose these equivalences to an equivalence between $\LN(\S_{\ord{0}})_{\ord{n}} = \li_0$ and $\S_{\ord{n}} = \li_n$.

An object of $\li_i$ is represented by a tuple
$$(C, \. x_{i+1}, x_{i+2}, \cdots, x_{n})$$
where $C$ is a $(i+2)$-complex and $x_{i+1}, \cdots, x_{n}$ are morphisms of $2$-complexes that satisfies
\begin{align*}
\mathrm{\bf v}_{i} C &= \dom x_{i+1}, \\
\cod x_{1} &= \dom x_{i+2}, \\
& \cdots \\
\cod x_{n-1} &= \dom x_{n}.
\end{align*}
Define the functor $\BG_{i} : \li_i \to \li_{i+1}$ as
$$(C, \. x_{i+1}, x_{i+2}, x_{i+3}, \cdots, x_{n}) \longmapsto (\G_i(C, x_{i+1}), x_{i+2} \epsilon''_{(C, x_{i+1}) 1}, x_{i+3}, \cdots, x_{n})$$
on objects and 
$$(a, \. f_{i+1}, f_{i+2}, f_{i+3}, \cdots, f_{n}) \longmapsto (\G_i(a, f_{i+1}), f_{i+2} \epsilon''_{(a, f_{i+1}) 1}, f_{i+3}, \cdots, f_{n})$$
on morphisms where the morphism $f_{i+2} \epsilon''_{(a, f_{i+1}) 1}$ is the composition

\[ \xymatrix@C=23pt@R=45pt{
\mathrm{\bf v}_{n+1} \G_n (C, x_{i+1}) \ar[d]_{\mathrm{\bf v}_{n+1} \G_n (a, f_{i+1})} \ar[rr]^(.65){\epsilon''_{(C, x_{i+1}) 1}} \ar@{.>}[drr]^(.45){\hat{\epsilon}''_{(a, f_{i+1}) 1}} \xtwocell[drr]{}\omit{^} && \bullet \ar[d]|{f_{i+2, 0}} \ar[rrr]^{x_{i+2}} \ar@{.>}[drrr]^(.45){\hat{f}_{i+2}} \xtwocell[drrr]{}\omit{^} &&& \bullet \ar[d]^{{f_{i+2, 1}}} \\
\mathrm{\bf v}_{n+1} \G_n (D, y_{i+1}) \ar[rr]_(.65){\epsilon''_{(D, y_{i+1}) 1}} && \bullet \ar[rrr]_{y_{i+2}} &&& \bullet
} \]

Define the functor $\BF_{i} : \li_{i+1} \to \li_i$ as
$$(C, \. x_{i+2}, x_{i+3}, \cdots, x_{n}) \longmapsto (\mathrm{\bf d}_{i+1} C, \tau, x_{i+2}, x_{i+3}, \cdots, x_{n})$$
where $\tau : \mathrm{\bf v}_i C \to \mathrm{\bf v}_{i+1} C$ and analogously for morphisms. It is easy to see that
$$\Spine_{\ord{n}} = \BF_0 \BF_1 \cdots \BF_{n-1}.$$
Define
$$\:\:\: \Filler_{\ord{n}} := \BG_{n-1} \BG_{n-2} \cdots \BG_0.$$

To prove the theorem, it is enough to prove that every $\BG_{i}$ is a weak inverse to $\BF_{i}$. By composing the equimorphisms, one gets the desired equivalence
\[ \xymatrix@C=25pt{
\LN(\S_{\ord{0}})_{\ord{n}} = \li_0 \ar@<.5ex>[r]^(.74){\BG_0} \ar@{<-}@<-.5ex>[r]_(.74){\BF_0} & \li_1 \ar@<.5ex>[r]^{\BG_1} \ar@{<-}@<-.5ex>[r]_{\BF_1} & \cdots \ar@<.5ex>[r]^(.38){\BG_{n-1}} \ar@{<-}@<-.5ex>[r]_(.38){\BF_{n-1}} & \li_n = \S_{\ord{n}}.}\]

Define natural transformations 
\[ \xymatrix@C=35pt{
\BG_i \BF_i \ar@<.5ex>[r]^(.5){\bar{\theta}} \ar@{<-}@<-.5ex>[r]_(.5){\bar{\eta}} & \Id_{\li_{i+1}}
} \]
as
\begin{align*}
\bar{\theta} & := (\theta, (\mathrm{\bf v}_{i+1} \theta, 1, 0), \id, \cdots, \id), \\
\bar{\eta} & := (\eta, (\mathrm{\bf v}_{i+1} \eta, 1, 0), \id, \cdots, \id)
\end{align*}
and lax natural transformations
\[ \xymatrix@C=35pt{
\BF_i \BG_i \ar@<.5ex>[r]^(.5){\bar{\epsilon}} \ar@{<-}@<-.5ex>[r]_(.5){\bar{\zeta}} & \Id_{\li_{i}}
} \]
as
\begin{align*}
\bar{\epsilon} & := (\epsilon', \epsilon'', (\epsilon''_1, 1, 0), \id, \cdots, \id), \\
\bar{\zeta} & := (\zeta', \zeta'', (\zeta''_1, 1, 0), \id, \cdots, \id).
\end{align*}

Now one has to check that the second component of $\bar{\theta}$ and $\bar{\eta}$ and the third component of $\bar{\epsilon}$ and $\bar{\zeta}$ are well-defined.
To verify that $(\epsilon''_1, 1, 0)$ and $(\zeta''_1, 1, 0)$ are well-defined, we see they map the object $(C, \. x_{i+1}, \cdots, x_{n})$ of $\li_{i}$ to 
\[\xymatrix@C=40pt@R=40pt{
\bullet \ar[d]_{\epsilon''_{(C, x_{i+1}) 1}} \ar[r]^{x_{i+2} \epsilon''_{(C, x_{i+1}) 1}} & \bullet \ar[d]^1 \\
\bullet \ar[r]_{x_{i+2}} & \bullet }
\xymatrix@C=20pt@R=15pt{ \\ & \mbox{ and } &}
\xymatrix@C=40pt@R=40pt{
\bullet \ar[d]_{\zeta''_{(C, x_{i+1}) 1}} \ar[r]^{x_{i+2}} & \bullet \ar@<-.6ex>[d]^1 \phantom{\: ,} \\
\bullet \ar[r]_{x_{i+2} \epsilon''_{(C, x_{i+1}) 1}} & \bullet \: ,}
\]
respectively. The squares clearly commute. The argument is analogous for morphisms.

Similarly for $(\mathrm{\bf v}_{i+1} \theta, 1, 0)$ and $(\mathrm{\bf v}_{i+1} \eta, 1, 0)$, observe that for objects $C$ of $\S_{\ord{i+1}}$
$$\epsilon''_{(\F_i C) 1} = \mathrm{\bf v}_{i+1} \theta_C$$
and analogously for morphisms (the left-hand side is natural).
The rest is analogous to the above.
\end{proof}

\begin{theorem}\label{adj_eq2}
The functors $\Spine_{\ord{n}}$ and $\Filler_{\ord{n}}$ are (part of) an adjoint equivalence.
\end{theorem}

\begin{proof}
Each tuple $\left( \BF_i, \BG_i, \bar{\epsilon}, \bar{\eta} \right)$ in Theorem \ref{th_filler} is an adjoint equivalence by verification similar to \ref{FG_adj}.
They are used to define $\Spine_{\ord{n}}$ and $\Filler_{\ord{n}}$ by composition of equimorphisms, but this composition respects composition of adjunctions as defined in \cite[IV.8]{Mac}.
\end{proof}

\nsection{The Augmented Case}

\begin{proposition}
The adjoint equivalence constituted by $\Spine_\ord{n}$ and $\Filler_\ord{n}$ can be extended to $n = -1$ and the simplicial functor $\Spine$ extended to its augmentation.
\end{proposition}

\begin{proof}
We need to define the functors $\Spine_\ord{-1}$ and $\Filler_\ord{-1}$, show that they are equimorphisms and that the diagram commutes
\[ \xymatrix{
\S_{\ord{0}} \ar[rr]^(.45){\Spine_{\ord{0}}} \ar[d]_{\mathrm{\bf d}_{0}} & &  \LN_+(\S_{\ord{0}})_{\ord{0}} \ar[d]^{\mathrm{\bf d}_{0}} \\
\S_{\ord{-1}} \ar[rr]^(.45){\Spine_{\ord{-1}}} & & \LN_+(\S_{\ord{0}})_{\ord{-1}}
} \]
where $\S = \Sigma \Com{\Delta}(\A)$. Recall that a d{\'e}calage of a simplicial object can be regarded as augmented by
$$\Sigma \Com{\Delta}(\A)_{\ord{-1}} := \Com{\Delta}(\A)_{\ord{0}} = \Com{1}(\A).$$
By definition
$$\LN_+(\Com{2}(\A))_{\ord{-1}} = \LFun( \ord{-1} , \: \Com{2}(\A)) \simeq \ord{0}.$$
Now define $\Filler_{\ord{-1}}$ to map the only object of $\ord{0}$ to the zero 1-complex in $\Com{1}(\A)$
\[ \xymatrix{
\Com{2}(\A) \ar[rr]^(.55){\Id} \ar[d]_{\mathrm{\bf d}_{1}} & & \Com{2}(\A) \ar[d] \\
\Com{1}(\A) \ar@<.5ex>[rr] & & \ar@<.5ex>[ll]^(.45){0} \phantom{.} \ord{0} .
} \]
The right going square commutes because there is only one functor $\Com{2}(\A) \to \ord{0}$.
It follows from Proposition \ref{com1} that the bottom part is an adjoint equivalence.
\end{proof}

\nsection{The Mapping Cone Functor}\label{cone}
Given a morphism $x : C \to X$ of chain complexes (i.e.\ 2-complexes), there is an associated chain complex $\Cone(x)$:
\[ \xymatrix@C=30pt{
\cdots \ar[r]^(.4){\mtrx{-d}{0}{x}{d}} & C^1 \oplus X^0 \ar[r]^{\mtrx{-d}{0}{x}{d}} & C^{2} \oplus X^1 \ar[r]^{\mtrx{-d}{0}{x}{d}} & C^{3} \oplus X^2 \ar[r]^{\mtrx{-d}{0}{x}{d}} & \cdots \: , 
} \]
defined in \cite[III.3.2]{GM}, called the mapping cone of $x$.
\begin{proposition}
If $x$ is regarded as an object of $\LArr(\Com{2}(\A))$ then
$$\Cone(x) = \mathrm{\bf d}_0 \Filler_\ord{1}(x),$$
where $\mathrm{\bf d}_0$ is the face functor of (the non-d{\'e}calaged) $\Com{\Delta}(\A)$.
\end{proposition}
\begin{proof}
Recall that $\Filler_\ord{1}(x) = \G_0(C, x)$ was the following 3-complex
\[ \xymatrix@C=20pt{
\cdots \ar[r] & C^0 \oplus X^0 \ar[r]^{\mtrx{d}{0}{0}{1}} & C^1 \oplus X^0 \ar[r]^(.55){\mtrx{-d}{0}{x}{d}} & C^2 \oplus X^1 \ar[r]^{\mtrx{-1}{0}{x}{d}} & C^2 \oplus X^2 \ar[r]^(.6){\mtrx{d}{0}{0}{1}} & \cdots
} \]
so we only need to verify that
$$\mtrx{d}{0}{0}{1} \mtrx{-1}{0}{x}{d} = \mtrx{-d}{0}{x}{d}.$$
\end{proof}
This result will be studied further in the next chapter.
At the moment, we extend the $\Cone$ to the whole sequicategory $\LArr(\Com{2}(\A))$.
\begin{definition}\label{mapping_cone}
Define the \textit{mapping cone functor}
$$\Cone : \LArr(\Com{2}(\A)) \to \Com{2}(\A)$$
as the composition 
\[ \xymatrix{
\Com{3}(\A) \ar[d]_{\mathrm{\bf d}_0} & & \ar[ll]_{\Filler_\ord{1}} \LArr(\Com{2}(\A)) \\
\Com{2}(\A) \ar[rr]^{\Spine_\ord{0} = \Id} & & \Com{2}(\A) .
} \]
where $\mathrm{\bf d}_0$ is a face functor of $\Com{\Delta}(\A)$.
\end{definition}

\section{Triangulations}\label{triangulations}
\renewcommand{\S}{{\mathcal S}}

\nsection{Semitriangulation}
Let $\C$ be a sesquicategory.
We begin by showing that a weak recalage of $\LN_+(\C)$ implies an ``semitriangulation'' on $\Ho(\C)$.

\begin{definition}\label{semitriangulation}
A nonempty category $\HH$, together with a collection of composable morphisms in $\HH$
$$\xymatrix@C=20pt{X \ar[r] & Y \ar[r] & Z}$$
called \textit{triangles} and a collection of objects in $\HH$ called \textit{vertices}, is said to be \textit{semitriangulated} if it satisfies the following axioms:
\begin{description}
 \item[S0] The collection of triangles and vertices are both closed under isomorphism, moreover, all vertices are isomorphic to each other.
 \item[S1] For every object $X$, there are triangles
 $$\xymatrix@C=20pt{X \ar[r]^{\id_X} & X \ar[r] & V}$$
 and
 $$\xymatrix@C=20pt{V \ar[r] & X \ar[r]^{\id_X} & X},$$
 which are unique up to isomorphism and where $V$ is a vertex.
 \item[S2] For any morphism $\xymatrix@C=20pt{X \ar[r]^f & Y}$, there exists a triangle
 $$\xymatrix@C=20pt{X \ar[r]^f & Y \ar[r] & Z}.$$
 \item[S4] For any commutative diagram 
 $$\xymatrix@C=20pt@R=20pt{A \ar[r] \ar[d]_{f_0} & B \ar[r] \ar[d]_{f_1} & C \\ D \ar[r] & E \ar[r] & F }$$
 with triangles at the rows, there exists an additional morphism $f_2$
 $$\xymatrix@C=20pt@R=20pt{A \ar[r] \ar[d]_{f_0} & B \ar[r] \ar[d]_{f_1} & C \ar@{.>}[d]^{f_2} \\ D \ar[r] & E \ar[r] & F }$$
 such that the diagram commutes. Moreover, if $f_0$ and $f_1$ are isomorphisms, then so is $f_2$.
 \item[S3] For any composable pair of morphisms $\xymatrix@C=20pt{X \ar[r]^f & Y \ar[r]^g & Z}$ there exist a commutative diagram
 \[\xymatrix@C=20pt@R=20pt{X \ar[r]^f \ar[dr] & Y \ar[rrr] \ar[d]^g & & & A \ar[ddll] \\
  & Z \ar[rd] \ar[dd] & & \\
  & & B \ar[dl] \\
  & C & & }\]
 where the four straight sub-diagrams are triangles.
\end{description}
\end{definition}

\begin{theorem}\label{catalmtr}
Let $\C$ be a sesquicategory. A recalage of $\LN_+(\C)$ induces a semitriangulation on $\Ho(\C)$.
\end{theorem}

\begin{proof}
Let $\S$ be a simplicial category such that $\LN_+(\C) = \Sigma \S$. In particular, $\C = \S_{\ord{1}}$. We shall construct a semitriangulated structure on $\Ho(\C)$.

First recall that a morphism in $\C$ is an equimorphism if and only if its equivalence class in $\Ho(\C)$ is an equimorphism, that is an isomorphism by \ref{isomorphism}.
Let the vertices of $\Ho(\C)$ be those objects in $\C$ equivalent to degeneracy of 0-simplices and
let the triangles of $\Ho(\C)$ be those diagrams in $\C$ equivalent to boundary of 2-simplices.
\vspace{-1em}
\[\xymatrix@C=50pt{
\S_{\ord{2}} \ar[r]|{\mathrm{\bf d}_1} \rtwocell^{\mathrm{\bf d}_2}_{\mathrm{\bf d}_0}<8>{\omit} \rtwocell\omit{<-2>^\tau} \rtwocell\omit{<2>^\tau} & \S_{\ord{1}} &  \ar[l]_{\mathrm{\bf s}_0} \S_{\ord{0}} \\
}\] $\.$
\vspace{-1em} \\
More explicitly, any object which is equivalent to $\mathrm{\bf s}_0 X$, for any $X$ in $\S_{\ord{0}}$, is a vertex and any diagram equivalent to one represented by
$$\xymatrix@C=30pt{\mathrm{\bf d}_2 Y \ar[r]^{\tau} & \mathrm{\bf d}_1 Y \ar[r]^{\tau} & \mathrm{\bf d}_0 Y},$$
for any $Y$ in $\S_{\ord{2}}$, is a triangle. Then {\bf S0} holds by construction.

To prove {\bf S1}, let $X$ be an object of $\S_{\ord{1}}$, then the triangle coming from $\mathrm{\bf s}_1 X$
$$\xymatrix@C=30pt{\mathrm{\bf d}_2 \mathrm{\bf s}_1 X \ar[r]^{\tau} & \mathrm{\bf d}_1 \mathrm{\bf s}_1 X \ar[r]^{\tau} & \mathrm{\bf d}_0 \mathrm{\bf s}_1 X}$$
is equal to
$$\xymatrix@C=30pt{X \ar[r]^{\id_X} & X \ar[r]^(.4){\tau} & \mathrm{\bf s}_0 \mathrm{\bf d}_0 X}$$
by the simplicial relations.

Axioms {\bf S2} and {\bf S3} follow from the equalities
\begin{align*}
\LN_+(\C)_\ord{1} &= \S_\ord{2}, \\
\LN_+(\C)_\ord{2} &= \S_\ord{3}
\end{align*}
respectively.
In this simplicial language, {\bf S2} and {\bf S3} describe the ``shadows'' of the 2- and 3-simplices of $\S$;   {\bf S3} is nothing but a description of a tetrahedron.

\newcommand{\lax}{\ar@{.>}[dr] \xtwocell[dr]{}\omit{_}}
\newcommand{\laz}{\ar@{.>}[ur] \xtwocell[ur]{}\omit{^}}

To prove {\bf S4}, assume we are given a commutative diagram in $\Ho(\C)$
$$\xymatrix@C=20pt@R=20pt{A \ar[r] \ar[d]_{f_0} & B \ar[r] \ar[d]_{f_1} & C \\ D \ar[r] & E \ar[r] & F }$$
with triangles at the rows. Choose representatives for this diagram in $\C$. By definition of triangles, there are objects $X, Y$ in $\S_{\ord{2}}$ 
\[\xymatrix@C=20pt@R=20pt{
\mathrm{\bf d}_2 X \ar[r]^{\tau} \ar@{=}[d] \lax & \mathrm{\bf d}_1 X \ar[r]^{\tau} \ar@{=}[d] \lax & \mathrm{\bf d}_0 X \ar@{=}[d] \\
A \ar[r] \ar[d] \lax & B \ar[r] \ar[d] & C \\
D \ar[r] \lax & E \ar[r] \lax & F \\
\mathrm{\bf d}_2 Y \ar[r]_{\tau} \ar@{=}[u] & \mathrm{\bf d}_1 Y \ar[r]_{\tau} \ar@{=}[u] & \mathrm{\bf d}_0 Y \ar@{=}[u] \\
}\]
such that the double arrows in the diagram are isomorphisms and the diagram commutes. Ignoring the rightmost column and using the equivalence between $\LArr(\C)$ and $\S_{\ord{2}}$, we find a corresponding $g : X \to Y$ in $\S_{\ord{2}}$. Taking $g$ back to $\C$ we get
\[\xymatrix@C=20pt@R=20pt{
A \ar[r] & B \ar[r] & C \\
\mathrm{\bf d}_2 X \ar[r] \ar@{=}[u] \ar[d]_(.45){\mathrm{\bf d}_2 g} \lax \laz & \mathrm{\bf d}_1 X \ar[r] \ar@{=}[u] \ar[d]|(.45){\mathrm{\bf d}_1 g \:} \lax \laz & \mathrm{\bf d}_0 X \ar@{=}[u] \ar[d]^(.45){\mathrm{\bf d}_0 g} \\
\mathrm{\bf d}_2 Y \ar[r] \ar@{=}[d] & \mathrm{\bf d}_1 Y \ar[r] \ar@{=}[d] & \mathrm{\bf d}_0 Y \ar@{=}[d] \\
D \ar[r] \laz & E \ar[r] \laz & F \\
}\]
The (equivalence class of the) composition $f_2 : C \to F$ is the sought morphism. To see that $f_2$ is an isomorphism when $f_0$ and $f_1$ are, note that the pair $(f_0, f_1)$ can be lifted to a equivalence in $\LArr(\C)$. Since $f_2$ is the value of a composition of equimorphism-preserving functors, it must be an isomorphism in $\Ho(\C)$.
\end{proof}

Recall that a pointed category is a category with zero object(s), that is objects which are simultaneously initial and terminal.
We always consider zero preserving functors between these categories.
Additive categories with additive functors are one example of pointed categories.

\begin{proposition}
Assume \. $\LN_+(\C)$ in Theorem \ref{catalmtr} is pointed. Then so is $\Ho(\C)$. Moreover, the vertices and zero objects of $\Ho(\C)$ coincide.
\end{proposition}

\begin{proof}
The category $\LN_+(\C)_{\ord{0}} = \C$ is pointed by assumption.
Any zero object of $\C$ remains a zero object of $\Ho(\C)$ because any groupoid with a single object has precisely one connected component.
Since $\LN_+(\C)_{\ord{-1}}$ has a single object, it is a zero object.
Let $\Sigma \S \simeq \LN_+(\C)$, then all objects of $\S_{\ord{0}}$ must be zero objects too. 
And since the functor $\mathrm{\bf s}_0 : \S_{\ord{0}} \longrightarrow \S_{\ord{1}}$ preserves zero objects, all vertices of $\Ho(\C)$ are zero objects.
\end{proof}

\begin{proposition}\label{stronger_axiom}
Let $\HH$ be a semitriangulated category.
Given a commutative diagram in $\HH$ represented by the solid arrows
\[\xymatrix@C=20pt@R=20pt{X \ar[r] \ar[dr] & Y \ar[rrr] \ar[d] & & & D \ar@{.>}[ddll] \\
 & Z \ar[rd] \ar[dd] & & \\
 & & E \ar@{.>}[dl] \\
 & F & & }\]
where the three straight solid lines are triangles, there exist two morphisms, represented by dotted arrows, that constitute a triangle and make the entire diagram commute.
\end{proposition}
\begin{proof}
Drop $D$, $E$, $F$ and apply {\bf S3}. By {\bf S4}, the objects $D$, $E$ and $F$ are isomorphic to $A$, $B$ and $C$, respectively. The dotted arrows are provided via the isomorphisms and the dotted line is a triangle by {\bf S0}.
\end{proof}

\nsection{Simplicial Triangulation}
We saw in the previous section that an recalage of $\LN_+(\C)$ provides a semitriangulated structure on  $\Ho(\C)$.

In this section we shall define a categorification of semitriangulation, called simplicial triangulation, and show that a simplicial triangulation on $\C$ is the same thing as an recalage of $\LN_+(\C)$.
Simplicial triangulation will also provide an enhancement of semitriangulation; it replaces some (weak) axioms with a ``mapping cone'' functor. 

\begin{definition}\label{simplicial_triangulation}
A \textit{simplicial triangulation} on a sesquicategory $\C$ is additional structure on the category given by a two functors, two lax natural transformations and a four relations.
The two functors are denoted $\mathrm{\bf d}_{-1}^\diamond$ and $\mathrm{\bf s}_{-1}^\diamond$:
\[ \xymatrix@C=35pt{
\LArr(\C) \ar[r]^(.6){\mathrm{\bf d}_{-1}^\diamond} & \C & \ar[l]_{\mathrm{\bf s}_{-1}^\diamond} \ord{0} 
} \]
called the \textit{mapping cone} and \textit{vertex}, respectively.
Identify
\begin{align*}
\ord{0} &= \LN_+(\C)_\ord{-1}, \\
\C &= \LN_+(\C)_{\ord{0}} , \\
\LArr(\C) &= \LN_+(\C)_{\ord{1}} .
\end{align*}
The two lax natural transformations are denoted $\tau^\mathrm{d}$ and $\tau^\mathrm{s}$:
\begin{align*}
\tau^\mathrm{d} &: \cod \longrightarrow \mathrm{\bf d}_{-1}^\diamond : \LArr(\C) \longrightarrow \C , \\
\tau^\mathrm{s} &: \mathrm{\bf s}_{-1}^\diamond \mathrm{\bf d}_0 \longrightarrow \id : \C \longrightarrow \C .
\end{align*}
Recall that $\cod = \mathrm{\bf d}_0 : \LArr(\C) \to \C$, so $\tau^\mathrm{d}$ can be thought of as the missing $\tau$
\[ \xymatrix@C=30pt{
\mathrm{\bf d}_1 \ar[r]^{\tau} & \mathrm{\bf d}_0 \ar[r]^{\tau^\mathrm{d}} & \mathrm{\bf d}_{-1} .
} \]
To explain $\tau^\mathrm{s}$, assume there would exist an $\mathrm{\bf s}_{-1} : \C \to \LArr(\C)$ that behaved according to the simplicial relations.
Then consider the diagram in $\Fun(\C, \C)$
\[ \xymatrix@C=30pt{
\mathrm{\bf d}_1 \mathrm{\bf s}_{-1} \ar[r]^{\tau} & \mathrm{\bf d}_0 \mathrm{\bf s}_{-1}
} \]
and use the simplicial relations $\mathrm{\bf d}_1 \mathrm{\bf s}_{-1} = \mathrm{\bf s}_{-1}^\diamond \mathrm{\bf d}_0$ and $\mathrm{\bf d}_0 \mathrm{\bf s}_{-1} = \id$, so $\tau^\mathrm{s}$ plays the role of $\tau$.
In fact, since we identify $\tau : \mathrm{\bf d}_1 \to \mathrm{\bf d}_0 : \LArr(\C) \to \C$ with the identity functor on $\LArr(\C) = \LN_+(\C)_{\ord{1}}$, we define
$$\tau^\mathrm{s} =: \mathrm{\bf s}_{-1} : \C \longrightarrow \LArr(\C).$$
Define
$$\mathrm{\bf d}_{-1} : \C \longrightarrow \ord{0}$$
as the only possible functor and 
$$\mathrm{\bf d}_{-1} : \LN_+(\C)_{\ord{2}} \longrightarrow \LArr(\C)$$
as given by the diagram in $\Fun(\LN_+(\C)_{\ord{2}}, \C)$
$$\xymatrix@C=30pt{\mathrm{\bf d}_{-1}^\diamond \mathrm{\bf d}_2 \ar[r]^{\tau = \mathrm{\bf d}_{-1}^\diamond \tau} & \mathrm{\bf d}_{-1}^\diamond \mathrm{\bf d}_1}.$$

The relations needed in the definition are those additional simplicial relations expressible in the recalage of the augmented sesquinerve $\LN_+(\C)$ involving $\mathrm{\bf d}_{-1}^\diamond$ and $\mathrm{\bf s}_{-1}^\diamond$, namely the following:
\begin{align*}
\mathrm{\bf d}_{-1}^\diamond \mathrm{\bf d}_0 \simeq \mathrm{\bf d}_{-1}^\diamond \mathrm{\bf d}_{-1} & : \LN_+(\C)_{\ord{2}} \longrightarrow \LN_+(\C)_{\ord{0}} & \mathrm{{\bf R1}} \\ 
\mathrm{\bf d}_{-1}^\diamond \mathrm{\bf d}_1 = \mathrm{\bf d}_0 \mathrm{\bf d}_{-1} & : \LN_+(\C)_{\ord{2}} \longrightarrow \LN_+(\C)_{\ord{0}} \\
\mathrm{\bf d}_{-1}^\diamond \mathrm{\bf d}_2 = \mathrm{\bf d}_1 \mathrm{\bf d}_{-1} & : \LN_+(\C)_{\ord{2}} \longrightarrow \LN_+(\C)_{\ord{0}} \\
\mathrm{\bf d}_{-1} \mathrm{\bf d}_0 = \mathrm{\bf d}_{-1} \mathrm{\bf d}_{-1}^\diamond & : \LN_+(\C)_{\ord{1}} \longrightarrow \LN_+(\C)_{\ord{-1}} \\
\mathrm{\bf d}_{-1} \mathrm{\bf d}_1 = \mathrm{\bf d}_0 \mathrm{\bf d}_{-1}^\diamond & : \LN_+(\C)_{\ord{1}} \longrightarrow \LN_+(\C)_{\ord{-1}} \\
\mathrm{\bf s}_{-1} \mathrm{\bf s}_{-1}^\diamond = \mathrm{\bf s}_0 \mathrm{\bf s}_{-1}^\diamond & : \LN_+(\C)_{\ord{-1}} \longrightarrow \LN_+(\C)_{\ord{1}} \\
\mathrm{\bf d}_{-1}^\diamond \mathrm{\bf s}_{-1} \simeq \mathrm{\bf d}_0 \mathrm{\bf s}_{-1} \simeq \id & : \LN_+(\C)_{\ord{0}} \longrightarrow \LN_+(\C)_{\ord{0}} & \mathrm{{\bf R2}} \\
\mathrm{\bf d}_{-1} \mathrm{\bf s}_{-1}^\diamond = \mathrm{\bf d}_0 \mathrm{\bf s}_{-1}^\diamond = \id & : \LN_+(\C)_{\ord{-1}} \longrightarrow \LN_+(\C)_{\ord{-1}} \\
\mathrm{\bf d}_{-1} \mathrm{\bf s}_1 = \mathrm{\bf s}_0 \mathrm{\bf d}_{-1}^\diamond & : \LN_+(\C)_{\ord{1}} \longrightarrow \LN_+(\C)_{\ord{1}} \\
\mathrm{\bf d}_{-1} \mathrm{\bf s}_0 \simeq \mathrm{\bf s}_{-1} \mathrm{\bf d}_{-1}^\diamond & : \LN_+(\C)_{\ord{1}} \longrightarrow \LN_+(\C)_{\ord{1}} & \mathrm{{\bf R3}} \\
\mathrm{\bf d}_1 \mathrm{\bf s}_{-1} = \mathrm{\bf s}_{-1}^\diamond \mathrm{\bf d}_0 & : \LN_+(\C)_{\ord{0}} \longrightarrow \LN_+(\C)_{\ord{0}} \\
\mathrm{\bf d}_{-1}^\diamond \mathrm{\bf s}_0 \simeq \mathrm{\bf s}_{-1}^\diamond \mathrm{\bf d}_{-1} & : \LN_+(\C)_{\ord{0}} \longrightarrow \LN_+(\C)_{\ord{0}} & \mathrm{{\bf R4}} 
\end{align*}
where $\simeq$ denotes the existence of an equimorphism. Note that all relation expressed with equality hold by construction and are therefore redundant. The remaining relations, labeled {\bf R1} -- {\bf R4}, are the ones we need in this definition. In fact, {\bf R4} follows from {\bf R3}, so it may be removed.
\end{definition}

\begin{example}\label{sim_tr_ch_com}
In this example, let $\C$ denote the sesquicategory of chain complexes $\Com{2}(\A)$ where $\A$ is an additive category.
Then $\C$ has a structure of simplicial triangulation given by
\begin{align*}
\mathrm{\bf d}_{-1}^\diamond &:= \Cone : \LArr(\C) \longrightarrow \C , \\
\mathrm{\bf s}_{-1}^\diamond &:= 0 : \ord{0} \longrightarrow \C , \\
\tau^\mathrm{d} &:= \cvec{0}{1} : \cod \longrightarrow \mathrm{\bf d}_{-1}^\diamond : \LArr(\C) \longrightarrow \C , \\
\tau^\mathrm{s} &:= 0 : \mathrm{\bf s}_{-1}^\diamond \mathrm{\bf d}_0 \longrightarrow \id : \C \longrightarrow \C ,
\end{align*}
where $\Cone$ is defined in \ref{cone}. The map $\tau^\mathrm{d}$ is equal to $\bar{\pi}$ in \cite[III.3.3]{GM}.
To verify that $\tau^\mathrm{d}_x : X \to \Cone x$ is a morphism in $\Com{2}(\A)$ for every object $x : C \to X$ in $\LArr(\C)$, we check that
$$\mtrx{-d}{0}{x}{d} \cvec{0}{1} = \cvec{0}{1} d.$$
To verify that $\tau^\mathrm{d}$ is a natural transformation $\cod \to \mathrm{\bf d}_{-1}^\diamond$, we check that for every morphism $f = (f_0, f_1, \hat{f}) : x \to y$
\[\xymatrix@C=35pt@R=35pt{
C \ar[d]_x \ar[r]^{f_0} \ar@{.>}[dr]^(.6){\hat{f}} \xtwocell[dr]{}\omit{_} & D \ar[d]^y \\
X \ar[r]_{f_1} & Y \\
}\]
$\cod f = f_1$, $\mathrm{\bf d}_{-1}^\diamond f = \Cone f = \mtrx{f_0}{0}{\hat{f}}{f_1}$ and the diagram commutes
\[\xymatrix@C=35pt@R=35pt{
X \ar[d]_{\cvec{0}{1}} \ar[r]^{f_1} & Y \ar[d]^{\cvec{0}{1}} \\
\Cone x \ar[r]_{\mtrx{f_0}{0}{\hat{f}}{f_1}} & \Cone y \. .\\
}\]
The existence of equimorphisms {\bf R1}--{\bf R3} can be verified in a similar way.
\end{example}

\begin{corollary}\label{triangles_threecomplex}
Triangles of $\LN_+(\Com{2}(\A))$ with the simplicial triangulation defined in \ref{sim_tr_ch_com} are equivalent to triangles of $\Sigma \Com{\Delta}(\A)$ (i.e.\ 3-complexes) with the simplicial triangulation given by the recalage $\Com{\Delta}(\A)$.
\end{corollary}
\begin{proof}
This is because the simplicial triangulation on $\LN_+(\Com{2}(\A))$ is coming from $\Com{\Delta}(\A)$. More precisely, because the mapping cone, defined in Definition \ref{mapping_cone} is a differential $\mathrm{\bf d}_0$ of $\Com{\Delta}(\A)$.
\end{proof}

We want to think about a simplicial triangulation on $\C$ as ``the same thing'' as a weak recalage of $\LN_+(\C)$. In the language of category theory, our desire can be formulated as follows.

\begin{conjecture}
Let $\C$ be a sesquicategory. Then the category of weak recalages of $\LN_+(\C)$ is equivalent to the category of simplicial triangulations on $\C$.
\end{conjecture}

But we have not defined these categories.
In fact, we are not really interested in these details.
Simplicial triangulation was introduced to bridge the gap between semitriangulation and weak recalage and a conceptual bridge is sufficient for this purpose.
The following theorem is intended to provide this.

\begin{theorem}\label{simtriang_is_recalage}
There is a one-to-one correspondence between weak recalages of $\LN_+(\C)$ (up to equivalence) and simplicial triangulations on a sesquicategory $\C$.
\end{theorem}

\begin{proof}
Given a weak simplicial sesquicategory $\S$ such that $\Sigma \S = \LN_+(\C)$, the required functors and lax natural transformations of the simplicial triangulation on $\C = \S_\ord{1}$ are defined as
\begin{align*}
\mathrm{\bf d}_{-1}^\diamond &:= \mathrm{\bf d}_{0}^\S : \S_\ord{2} \longrightarrow \S_\ord{1} , \\
\mathrm{\bf s}_{-1}^\diamond &:= \mathrm{\bf s}_{0}^\S : \S_\ord{0} \longrightarrow \S_\ord{1} , \\
\tau^\mathrm{d} &:= \tau : \mathrm{\bf d}_1^\S \longrightarrow \mathrm{\bf d}_0^\S : S_\ord{2} \longrightarrow \S_\ord{1} , \\
\tau^\mathrm{s} &:= \tau : \mathrm{\bf d}_2^\S {\bf s}_0^\S \longrightarrow \mathrm{\bf d}_1^\S {\bf s}_0^\S : \S_\ord{1} \longrightarrow \S_\ord{2} .
\end{align*}
They satisfy the required relations because $\S$ was a weak simplicial category to begin with.

For the reverse direction, given a simplicially triangulated sesquicategory $\C$, define 
$$\S_\ord{i} := \LN_+(\C)_\ord{i-1}, \; i \geq 0.$$
Define the face and degeneracy functors $\mathrm{\bf d}_{i}^\S$ and $\mathrm{\bf s}_{i}^\S$ of $\S$ as
\begin{align*}
\mathrm{\bf d}_{i}^\S &:= \mathrm{\bf d}_{i-1}^{\LN_+(\C)}, \\
\mathrm{\bf s}_{i}^\S &:= \mathrm{\bf s}_{i-1}^{\LN_+(\C)}
\end{align*}
for $i \neq 0$ where $\mathrm{\bf d}_{i-1}^{\LN_+(\C)}$ and $\mathrm{\bf s}_{i-1}^{\LN_+(\C)}$ are the face and degeneracy functors for $\LN(\C)_+$.
This guaranties that $\Sigma \S = \LN_+(\C)$.
It remains to define the $\mathrm{\bf d}_0^\S$ and $\mathrm{\bf s}_{0}^\S$ functors.
For the case
\[ \xymatrix@C=30pt{
\S_{\ord{0}} \ar@<.6ex>[r]^{\mathrm{\bf d}_0^\S} \ar@<-.6ex>@{<-}[r]_{\mathrm{\bf s}_0^\S} & \S_{\ord{-1}}
} \]
let $\mathrm{\bf d}_0^\S$ be the unique constant functor and $\mathrm{\bf s}_0^\S := \mathrm{\bf s}_{-1}^\diamond$.
For
\[ \xymatrix@C=30pt{
\S_{\ord{1}} \ar@<.6ex>[r]^{\mathrm{\bf d}_0^\S} \ar@<-.6ex>@{<-}[r]_{\mathrm{\bf s}_0^\S} & \S_{\ord{0}}
} \]
let $\mathrm{\bf d}_0^\S := \mathrm{\bf d}_{-1}^\diamond$ and $\mathrm{\bf s}_0^\S := \mathrm{\bf s}_{-1} = \tau^{\mathrm s}$. For
\[ \xymatrix@C=30pt{
\S_{\ord{i+1}} \ar[r]^(.55){\mathrm{\bf d}_0^\S} & \S_{\ord{i}}
} \]
where $i \geq 1$, write $\S_{\ord{i}} = \LN_+(\C)_\ord{i-1}$ and define $\mathrm{\bf d}_0^\S$ as given by
$$\xymatrix{\mathrm{\bf e}_0 \ar[r]^\tau & \mathrm{\bf e}_1 \ar[r]^\tau & \dots \ar[r]^(.45)\tau & \mathrm{\bf e}_{i-1}}$$
where
$$\mathrm{\bf e}_j := \mathrm{\bf d}_{-1}^\diamond \mathrm{\bf d}_1 \cdots \mathrm{\bf d}_{j} \mathrm{\bf d}_{j+2} \cdots \mathrm{\bf d}_{i} : \S_{\ord{i+1}} \longrightarrow \C,$$
and the functors $\mathrm{\bf d}_j = \mathrm{\bf d}_j^{\LN_+(\C)}$.
For
\[ \xymatrix@C=30pt{
\S_{\ord{i+1}} & \ar[l]_(.4){\mathrm{\bf s}_0^\S} \S_{\ord{i}}
} \]
define $\mathrm{\bf s}_0^\S$ as given by
$$\xymatrix{\mathrm{\bf g}_0 \ar[r]^\tau & \mathrm{\bf g}_1 \ar[r]^\tau & \dots \ar[r]^(.45)\tau & \mathrm{\bf g}_i}$$
where
$$\mathrm{\bf g}_j := \mathrm{\bf d}_0 \cdots \mathrm{\bf d}_{j-2} \mathrm{\bf d}_{j} \cdots \mathrm{\bf d}_{i-1} : \S_{\ord{i}} \longrightarrow \C$$
for $j \geq 1$ and
$$\mathrm{\bf g}_0 := \mathrm{\bf s}_{-1}^\diamond \mathrm{\bf d}_0 \mathrm{\bf d}_1 \cdots \mathrm{\bf d}_{i-1} : \S_{\ord{i}} \longrightarrow \C$$
and $\tau : \mathrm{\bf g}_0 \to \mathrm{\bf g}_1$ is equal to $\tau^{\mathrm s} \. \mathrm{\bf d}_1 \cdots \mathrm{\bf d}_{i-1}$. The corresponding lax natural transformations are defined in the obvious way using $\tau^{\mathrm d}$ and $\tau^{\mathrm s}$.

This makes $\S$ into a simplicial sesquicategory, because all additional simplicial relations are satisfied by definition of simplicial triangulation. These are essentially $\mathrm{{\bf R1}} - \mathrm{{\bf R3}}$.
\end{proof}

\begin{theorem}\label{catagorifies}
Simplicial triangulation categorifies semitriangulation.
\end{theorem}

\begin{proof}
It follows from Theorems \ref{catalmtr} and \ref{simtriang_is_recalage} that a simplicial triangulation on $\C$ implies a semitriangulation on $\Ho(\C)$.
To prove this in accordance to Definition \ref{def_decat} we need to check that a semitriangulation is precisely the structure that remains after the decategorification $\C \mapsto \Ho(\C)$.
By this we mean that the data and relations in the definition of a simplicial category can be described in terms of categorification of the corresponding structure of semitriangulation. 

The functor $\mathrm{\bf s}_{-1}^\diamond : \ord{0} \to \C$ determines a vertex and consequently the isomorphism class of vertices.

Triangles, and axioms {\bf S2} and {\bf S4} are categorified as the functor 
$$\mathrm{\bf d}_{-1}^\diamond : \LArr(\C) \to \C$$
and lax natural transformation $\tau^{\mathrm d}$. That is, any triangle is isomorphic to one of the form 
\[ \xymatrix@C=30pt{
\mathrm{\bf d}_1 f \ar[r]^{\tau} & \mathrm{\bf d}_0 f \ar[r]^{\tau^\mathrm{d}_f} & \mathrm{\bf d}_{-1}^\diamond f 
} \]
or equivalently
\[ \xymatrix@C=30pt{
\dom f \ar[r]^{f} & \cod f \ar[r]^{\tau^\mathrm{d}_f} & \mathrm{\bf d}_{-1}^\diamond f 
} \]
where $f$ is an object of $\LArr(\C)$. Axiom {\bf S4} is an attempt to say that the object $Z$ given in {\bf S2} is the result of a decategorified functor (the mapping cone).

The natural transformation $\tau^{\mathrm s}$ and relations {\bf R2} and {\bf R4} categorifies the second triangle in {\bf S1}:
$$\xymatrix@C=30pt{V \ar[r]^{\tau^{\mathrm s}_X} & X \ar[r]^{\id_X} & X}.$$
Relation {\bf R3} can be seen as the uniqueness requirement of {\bf S1}.

Finally, relation {\bf R1} follows from axiom {\bf S3}.
\end{proof}

Informally, a semitriangulated structure on $\Ho(\C)$ is an attempt to capture the shadows casted by an recalage of $\LN_+(\C)$ onto the sesquicategory $\C$. This additional structure on $\LN_+(\C)$ behaves like an extra dimension:
\begin{itemize}
\item The $(-1)$-simplex takes on the role of a zero object.
\item The $0$-simplices, the objects of $\C$, become edges with a source and a target, this is reflected in the alternating sign of the Euler characteristic.
\item The $1$-simplices, the morphisms of $\C$, become triangles, where the newly discovered third edge behaves as a mapping cone in a functorial manner.
\item The $2$-simplices, composable pairs of morphisms, become tetrahedrons, with the four faces corresponding to four triangles (Verdier's axiom).
\end{itemize}

\nsection{Ordinary Triangulation}
\newcommand{\Tr}{\operatorname{T}}

For comparison, we include a definition of triangulated categories.

\begin{definition}\label{triangulation}
A \textit{triangulation} on an additive category $\HH$ is an additive self-equivalence $\Tr : \HH \to \HH$ together with a collection of diagrams of the form
$$\xymatrix@C=22pt{X \ar[r]^f & Y \ar[r]^g & Z \ar[r]^(.4)h & \Tr X}$$
called \textit{triangles} and denoted $(f, g, h)$, such that the following axioms hold.

\begin{description}
 \item[T0] The collection of triangles are closed under isomorphism.
 \item[T1] For any object $X$ in $\HH$, there is a triangle
 $$\xymatrix@C=22pt{X \ar[r]^1 & X \ar[r] & 0 \ar[r] & \Tr X .}$$
 \item[T2] For any morphism $f : X \to Y$ in $\HH$, there is a triangle $(f, g, h)$
 $$\xymatrix@C=22pt{X \ar[r]^f & Y \ar[r]^g & Z \ar[r]^(.4)h & \Tr X .}$$
 \item[T4] If the rows are triangles and the left square commutes in the following diagram, then there is a morphism $k$ that makes the remaining squares commute.
 \[\xymatrix@C=22pt{
 X \ar[r]^f \ar[d]_i & Y \ar[r]^g \ar[d]^j & Z \ar[r]^(.4)h \ar@{-->}[d]^k & \Tr X \ar[d]^{\Tr i} \\
 X' \ar[r]^{f'} & Y' \ar[r]^{g'} & Z' \ar[r]^(.42){h'} & \Tr X'
 }\]
 \item[T3] Consider the following diagram.
 \[\xymatrix@C=22pt@R=22pt{
 X \ar[dr]_f \rrlowertwocell<8>^{h}{\omit} && Z \ar@{-->}[dr]^{h'} \rrlowertwocell<8>^{g'}{\omit} && W \ar[dr]^{g''} \rrlowertwocell<8>^{j''}{\omit} && \Tr U \\
 & Y \ar[ur]^g \ar[dr]_{f'} && V \ar@{-->}[ur]^{j'} \ar@{-->}[dr]^{h''} && \Tr Y \ar[ur]_{\Tr f'} & \\
 && U \ar@{-->}[ur]^{j} \rrlowertwocell<-8>_{f''}{\omit} && \Tr X \ar[ur]_{\Tr f} && \\
 }\]
 Assume that $h = gf$, $j'' = (\Tr f') g''$, and $(f, f', f'')$ are $(g, g', g'')$ are triangles.
 If $h'$ and $h''$ are given such that $(h, h', h'')$ is a triangle, then there are morphisms $j$ and $j'$ such that the diagram commutes and $(j, j',  j'')$ is a triangle.
 \item[T5] Given a triangle
 $$\xymatrix@C=22pt{X \ar[r]^f & Y \ar[r]^g & Z \ar[r]^(.4)h & \Tr X}$$
 then the following is a triangle
 $$\xymatrix@C=22pt{Y \ar[r]^g & Z \ar[r]^(.4)h & \Tr X \ar[r]^(.4){-\Tr f} & \Tr Y.}$$
\end{description}

\end{definition}

Apart from {\bf T4}, this definition is from \cite{May}. Traditionally, {\bf T4} has been part of the definition, but in \cite{May}, May proves that {\bf T4} follows from the other axioms, hence it is not included in the definition. However, with our simplicial bias the proper course of action is to weaken {\bf T3} to {\bf S3} rather than removing {\bf T4}. This is because {\bf S3} is the higher analog of {\bf S2}. The stronger version of {\bf S3} then follows by Proposition \ref{stronger_axiom}.

Axiom {\bf T3} is known as Verdier's axiom or the octahedral axiom.
There is a surprising number of ways to draw the associated diagram of this axiom.
The term octahedron is due to the fact that one of them has named shape.
If a Platonic solid must name this axiom, then the tetrahedron is a better choice since it refers to its intrinsic properties.
Triangles are at least called triangles, but perhaps for the wrong reason.
According to \cite[Remark 1.3.5]{We}, the term stems from the shape they are often depicted in.
But in those diagrams the vertices are the objects of $\HH$ which, from our point of view, are the edges of a triangle.
Many authors even call any triangle-shaped diagram a triangle and call the true ones for distinguished triangles.

We now compare triangilation defined in \ref{triangulation} with semitriangulation defined in  \ref{semitriangulation}.
The following axioms need to be added to semitriangulation in order to make the two definitions equivalent:
\begin{description}
 \item[A0] The category $\HH$ is additive.
 \item[A1] There is a functor $\Tr : \HH \to \HH$.
 \item[A2] The functor $\Tr$ is an additive self-equivalence.
 \item[A3] Let
 $\xymatrix@C=22pt{X \ar[r]^f & Y \ar[r]^(.4)g & Z}$
 be a triangle.
 Then there is a morphism $h$ such that
 $$\xymatrix@C=22pt{Y \ar[r]^g & Z \ar[r]^(.4)h & \Tr X}$$
 is a triangle
 \item[A4] and
 $$\xymatrix@C=22pt{Z \ar[r]^(.4)h & \Tr X \ar[r]^(.42){-\Tr f} & \Tr Y}$$
 is a triangle.
 \item[A5] Given a morphism of triangles
 \[\xymatrix@C=22pt{
 X \ar[r]^f \ar[d]_i & Y \ar[r]^g \ar[d]^j & Z \ar[d]^k \\
 X' \ar[r]^{f'} & Y' \ar[r]^{g'} & Z'
 }\]
 then
 \[\xymatrix@C=22pt{
 Y \ar[r]^g \ar[d]^j & Z \ar[r]^(.4)h \ar[d]^k & \Tr X \ar[d]^{\Tr i} \\
 Y' \ar[r]^{g'} & Z' \ar[r]^(.42){h'} & \Tr X'
 }\]
 is a morphism of triangle.
\end{description}
Axiom {\bf A1} is the only one stipulating extra structure, but contrary to what it may seem, this axiom is the least crucial one. In fact we may augment semitriangulation with he following three axioms without affecting Theorem \ref{catagorifies} which says that simplicial triangulation categorifies semitriangulation.
\begin{description}
 \item[A1] There is functor $\Tr : \HH \to \HH$.
 \item[A3'] Let
  $\xymatrix@C=22pt{X \ar[r]^\id & X \ar[r]^v & V}$
  be a triangle.
  Then there is a triangle
  $$\xymatrix@C=22pt{X \ar[r]^v & V \ar[r] & \Tr X.}$$
 \item[A5'] Given a morphism of triangles
 \[\xymatrix@C=22pt{
 X \ar[r]^\id \ar[d]_i & X \ar[r]^v \ar[d]^i & V \ar[d]^k \\
 X' \ar[r]^\id & X' \ar[r]^{v'} & V'
 }\]
 then
 \[\xymatrix@C=22pt{
 X \ar[r]^v \ar[d]^i & V \ar[r]^(.4)h \ar[d]^k & \Tr X \ar[d]^{\Tr i} \\
 X' \ar[r]^{v'} & V' \ar[r]^(.42){h'} & \Tr X'
 }\]
 is a morphism of triangles.
\end{description}
To see this, consider 
$$\tau^\mathrm{d} : \cod \longrightarrow \mathrm{\bf d}_{-1}^\diamond : \LArr(\C) \longrightarrow \C,$$
from the definition of simplicial triangulation (Definition \ref{simplicial_triangulation}), as a functor \\ $\LArr(\C) \to \LArr(\C)$ and let $\Tr$ correspond to the decategorification of
$$\mathrm{\bf d}_{-1}^\diamond \tau^\mathrm{d} \mathrm{\bf s}_0 : \C \longrightarrow \C,$$
that is $\Tr := \Ho(\mathrm{\bf d}_{-1}^\diamond \tau^\mathrm{d} \mathrm{\bf s}_0)$.
In this way, we have access to (a candidate) of $\Tr$ in semitriangulated categories which may, or may not, satisfy {\bf A2} -- {\bf A4}.

\section{Groups}\label{groups}
In the previous chapter we saw that the categorification of semitriangulation is simplicial triangulation, which is the same thing as a recalage of an augmented lax nerve.
In this chapter we show that the decategorification of a semitriangulation is a group structure.
We then show that a group structure on a set $M$ is the same thing as a recalage of the augmented nerve of $M$.

\nsection{Division}

\begin{definition}\label{group}
A \textit{group} is a set $M$ together two functions 
\[ \xymatrix@C=35pt{
M \times M \ar[r]^(.6){\mathrm{div}} & M & \ar[l]_{e} \{ \emptyset \}
} \]
satisfying the following axioms for all $x, y, z \in M$
\begin{description}
 \item[L1] $x/x = e,$
 \item[L2] $x/e = x,$
 \item[L3] $(z/x)/(y/x) = z/y.$
\end{description}
where $/$ is infix for $\mathrm{div}$ and $e$ is shorthand for $e(\emptyset)$.
\end{definition}

This definition appears in literature (see e.g. \cite[p. 6]{Ha}) but with the additional axiom
\begin{description}
 \item[L4] $e/(y/x) = x/y.$
\end{description}
However, {\bf L4} is redundant because it follows from {\bf L1} and {\bf L3}.

\begin{proposition}
The definition of group in Definition \ref{group} is equivalent to the classical one:
A group is a set $M$, together with a function $\_ * \! \_ : M \times M \to M$, a function $\_^{-1} : M \to M$ and an element $e \in M$ that satisfy for all $x, y, z \in M$:
\begin{description}
 \item[G1] $(x * y) * z = x * (y * z).$
 \item[G2] $e * x = x,$
 \item[G3] $x^{-1} * x = e,$
\end{description}
\end{proposition}
\begin{proof}
To prove that a {\bf G}-group is an {\bf L}-group, define
$$x/y := x * y^{-1}$$
and check {\bf L3}:
$$(z/x)/(y/x) = z * x^{-1} * (y * x^{-1})^{-1} = z* x^{-1} * x * y^{-1} = z/y.$$

To prove that an {\bf L}-group is a {\bf G}-group, define
\begin{align*}
x^{-1}  &:= e / x, \\
x * y &:= x / y^{-1}.
\end{align*}
Then {\bf G3} follows from {\bf L1}.
We prove {\bf G2}:
$$e * x = e / (e / x) = (x / x) / (e / x) = x / e = x$$
and {\bf G1}:
\begin{align*}
(x * y) * z &= (x * y) / (e / z) \\
&= (x / (e / y)) / (e / z) \\
&= (x / (e / y)) / ((e / z) / e) \\
&= (x / (e / y)) / (((e / z) / y) / (e / y)) \\
&= x / ((e / z) / y) \\
&= x / (((e / z) / (e / z)) / (y / (e / z))) \\
&= x / (e / (y / (e / z))) \\
&= x * (y / (e / z)) \\
&= x * (y * z) .
\end{align*}
\end{proof}

\nsection{Recalage of the Nerve of a Set}

\begin{theorem}\label{alm_tr_cat_group}
The structure of semitriangulation on a category categorifies a group structure on a set.
\end{theorem}
\begin{proof}
Decategorifying semitriangulation defined in Definition \ref{semitriangulation}, we get the following:

A nonempty set $M$, together with a collection of ``triangles''
$$T \subset M \times M \times M$$
and a collection of ``vertices'' in $M$ that satisfies the following axioms:
\begin{description}
 \item[$\Ho$(S0)] All vertices are equal (and they exist by {\bf S1}). This is our $e$.
 \item[$\Ho$(S2)] For any pair $(x, y) \in M \times M$ there is a $z \in M$ such that $(x, y, z) \in T$.
 \item[$\Ho$(S4)] Let $(a, b, c), (d, e, f) \in T$, if $(a, b) = (d, e)$ then $c = f$.
 Axioms $\Ho${\bf (S2)} and $\Ho${\bf (S4)} say precisely that $T$ is (the graph of) a function $M \times M \to M$. This function is our $\mathrm{div}$, i.e.\ $T = \{(x, y, y/x) \: | \: x, y \in M \}$.
 \item[$\Ho$(S1)] For any $x \in M$, we have $(x, x, e) \in T$ and $(e, x, x) \in T$, this is axioms {\bf L1} and {\bf L2}.
 \item[$\Ho$(S3)] For any $x, y, z \in M$, there exist $a, b, c \in M$ such that
\begin{align*}
(x, y, a) & \in T, \\
(x, z, b) & \in T, \\
(y, z, c) & \in T, \\
(a, b, c) & \in T.
\end{align*}
This is precisely axiom {\bf L3}.

Since all of the group structure as given in Definition \ref{group} is accounted for, this is a categorification which lives up to Definition \ref{def_decat}.
\end{description}
\end{proof}

\begin{corollary}\label{simtriang_cat_cat_group}
The structure of simplicial triangulation on a sesquicategory is a double categorification of the structure of a group on a set.
\end{corollary}
\begin{proof}
In the spirit of Definition \ref{def_decat}, this should follow from Theorem \ref{catagorifies} and \ref{alm_tr_cat_group}.
To be on safe side, we still go trough the details.
Let 
\[ \xymatrix@C=35pt{
M \times M \ar[r]^(.6){\mathrm{div}} & M & \ar[l]_{e} \{ \emptyset \}
} \]
be given by $\Ho^2 := \Ho \circ \Ho$ applied to
\[ \xymatrix@C=35pt{
\LArr(\C) \ar[r]^(.6){\mathrm{\bf d}_{-1}^\diamond} & \C & \ar[l]_{\mathrm{\bf s}_{-1}^\diamond} \ord{0}.
} \]
And for the relations
\begin{align*}
\Ho^2(\mathrm{{\bf R1}}) &= \mathrm{{\bf L3}}, \\
\Ho^2(\mathrm{{\bf R2}}) &= \mathrm{{\bf L2}}, \\
\Ho^2(\mathrm{{\bf R4}}) &= \mathrm{{\bf L1}}.
\end{align*}
In Definition \ref{simplicial_triangulation}, the relation $\mathrm{{\bf R4}}$ follows from $\mathrm{{\bf R3}}$.
But for the purpose of defining groups, it is more economical to keep $\Ho^2(\mathrm{{\bf R4}})$ and observe that $\Ho^2(\mathrm{{\bf R3}})$ follows from $\Ho^2(\mathrm{{\bf R4}})$, that is
$$(x/x, y/x) = (0, y/x)$$
for all $(x, y) \in M \times M$.
\end{proof}

We have shown that the double decategorification of simplicial triangulation is a group.
But a simplicial triangulation on $\C$ is equivalent to a weak recalage of the lax nerve of sesquicategory $\C$.
If decategorification in Definition \ref{def_decat} is well-behaved, it should follow that the group in question is in fact a recalage of the augmented nerve of the set $\Ho^2(\C) = \Ho(\Ho(\C))$.
This makes us suspect the following.

\begin{theorem}\label{recalage_set}
A group structure on a set $M$ uniquely determines a recalage of $\NN_+(M)$ and vice versa.
\end{theorem}
\begin{proof}
Assume we are given a simplicial set $\S$ such that $\Sigma \S = \NN_+(M)$.
Then $\S_\ord{n} = \NN_+(M)_\ord{n-1} = M^{\times n}$ is the $n$:th Cartesian power of $M$. Define
\begin{align*}
\mathrm{div} &:= \mathrm{\bf d}_{0} : \S_\ord{2} \longrightarrow \S_\ord{1} , \\
e &:= \mathrm{\bf s}_{0} : \S_\ord{0} \longrightarrow \S_\ord{1}.
\end{align*}
Since $\mathrm{\bf s}_{1} : \S_\ord{1} \to \S_\ord{2}$ is the diagonal map $M \to M \times M$, the relation 
$$\mathrm{\bf d}_{0} \mathrm{\bf s}_{1} = \mathrm{\bf s}_{0} \mathrm{\bf d}_{0} : \S_\ord{1} \longrightarrow \S_\ord{1}$$
means that for $x \in M$, $\mathrm{\bf d}_{0} \mathrm{\bf s}_{1} x = \mathrm{div}(x, x) = e$ which is Axiom {\bf L1}.
The relations
\begin{align*}
\mathrm{\bf d}_{1} \mathrm{\bf s}_{0} = \id &: \S_\ord{1} \longrightarrow \S_\ord{1} \\
\mathrm{\bf d}_{2} \mathrm{\bf s}_{0} = \mathrm{\bf s}_{0} \mathrm{\bf d}_{1} &: \S_\ord{1} \longrightarrow \S_\ord{1}
\end{align*}
tell us that for $x \in M$, $\mathrm{\bf s}_{0} x = (e, x)$ because $\mathrm{\bf d}_{1}$ and $\mathrm{\bf d}_{2}$ are the two projections $M \times M \to M$.
The relation
$$\mathrm{\bf d}_{0} \mathrm{\bf s}_{0} = \id : \S_\ord{1} \longrightarrow \S_\ord{1}$$
then means that for $x \in M$, $\mathrm{\bf d}_{0} \mathrm{\bf s}_{0} x = \mathrm{div}(e, x) = x$ which is Axiom {\bf L2}. Axiom {\bf L3} follows in similar way from the relation
$$\mathrm{\bf d}_0 \mathrm{\bf d}_1 = \mathrm{\bf d}_0 \mathrm{\bf d}_0 : \S_{\ord{3}} \longrightarrow \S_{\ord{1}}.$$
This proves that a recalage of $\NN_+(M)$ endows $M$ with a group structure.

To prove the opposite direction, let $(M, \mathrm{div}, e)$ be a group as defined in Definition \ref{group} and define a recalage $\S$ of $\NN_+(M)$ just like in Theorem \ref{simtriang_is_recalage}.
In particular
\begin{align*}
\S_\ord{i} &:= \NN_+(M)_\ord{i-1} = M^{\times i} & i &\geq 0 ,\\
\mathrm{\bf d}_{i}^\S &:= \mathrm{\bf d}_{i-1}^{\NN_+(M)} & i &\geq 1 , \\
\mathrm{\bf s}_{i}^\S &:= \mathrm{\bf s}_{i-1}^{\NN_+(M)} & i &\geq 1 
\end{align*}
and
\[ \xymatrix@C=40pt{
\S_\ord{2} \ar[r]^{\mathrm{\bf d}_0^\S := \mathrm{div}} & \S_\ord{1} & \ar[l]_{\mathrm{\bf s}_0^\S := e} \S_\ord{0} .
} \]

However, the reason $\S$ satisfies the simplicial relations differs slightly compared to Theorem \ref{simtriang_is_recalage}, this is because $\Ho^2(\mathrm{{\bf R3}})$ is redundant here as explained in Corollary \ref{simtriang_cat_cat_group}.
\end{proof}

\begin{definition}
Let $\mathrm{\bf Rec_1}$ denote the full subcategory of the category of simplicial sets, given by those simplicial sets which are recalage of nerves.
\end{definition}

In our 2-categorical terminology, these objects are functors $F : \Delta \to \mathrm{R} \circ \Set$ such that $\Sigma F = \NN_+(F_\ord{1})$, where $\mathrm{R}$ denotes the right adjoint to $\ob$ as described in Definition \ref{object}.

\begin{theorem}
The category $\mathrm{\bf Rec_1}$ is isomorphic to $\Grp$; the category of small groups.
\end{theorem}
\begin{proof}
The bijection on objects is given by Theorem \ref{recalage_set}.
Essentially, a simplicial set is turned into a group (see Definition \ref{group}) by chopping of degrees higher than two.
Hence every natural transformation preserves the group structure and gives a group homomorphism. This clearly defines a functor
$$\operatorname{chop} : \mathrm{\bf Rec_1} \longrightarrow \Grp.$$

To see that $\operatorname{chop}$ is faithful, note that every morphism $\eta : F \to G$ of $\mathrm{\bf Rec_1}$ is uniquely determined by $\eta_\ord{1}$ because the vertex maps (see Definition \ref{vertex}) of $\Sigma F$ are precisely the projections.

To see that $\operatorname{chop}$ is full, let $h : \operatorname{chop} F \to \operatorname{chop} G$ be a morphism in $\Grp$.
Then $h$ has a unique extension to a simplicial map $h' : \Sigma F \to \Sigma G$ for the same reason as above.
It remains to prove that $h'$ commutes with the remaining generators, that means
$$\mathrm{\bf d}_0^G h' = h' \mathrm{\bf d}_0^F$$
and analogously for $\mathrm{\bf s}_0$.
But $\mathrm{\bf d}_0$ and $\mathrm{\bf s}_0$ are defined in terms of projections, division and insertion of the neutral element (see the proof of Theorem \ref{simtriang_is_recalage}) which are respected by the group homomorphism $h$.
\end{proof}

\nsection{Nerve of a Group}
In this section we shall explain the difference between the recalage description of a group and the more familiar nerve description.
Since the recalage description is a decategorification of some of our other constructions, we will gain some insight into its simplicial structure.

Let 
$$G = (M, \mathrm{div}, e)$$
be a group as in Definition \ref{group} and $\S$ be the recalage of the nerve of $M$ that corresponds to $G$ as given by Theorem \ref{recalage_set}.
In other words, $\S$ is a simplicial set such that $\Sigma \S = \NN_+(M)$ and the diagram
\[ \xymatrix@C=40pt{
\S_\ord{2} \ar[r]^{\mathrm{\bf d}_0^\S} & \S_\ord{1} & \ar[l]_{\mathrm{\bf s}_0^\S} \S_\ord{0}
} \]
is equal to
\[ \xymatrix@C=35pt{
M \times M \ar[r]^(.6){\mathrm{div}} & M & \ar[l]_{e} \{ \emptyset \} .
} \]

On the other hand, we may regard $G$ as a groupoid with a single object and $M$ as its morphism set with composition given by the group structure.
Denote the simplicial set $\ob(\NN(G))$, the nerve of the groupoid $G$, by $\mathrm{N}(G)$.

Then both $\S$ and $\mathrm{N}(G)$ are simplicial sets that fully describe the group $G$.
But they do that in different ways.
For example, most of the generators of the simplicial maps of $\S$ (those of $\Sigma \S$ to be specific) has no knowledge of the group structure while the rest have explicit knowledge of division.
In contrast, the simplicial maps of $\mathrm{N}(G)$ have explicit knowledge of composition inside $G$, but only implicit knowledge of division; $G$ could have been a monoid (a category with a single object) in the later construction.
The easiest way to understand the difference is trough a picture.

Let $(x, y, z)$ be an ordered triple of elements of the group $G$. The triple then corresponds to a unique 3-simplex (tetrahedron) in both $\S$ and $\mathrm{N}(G)$:
\newcommand{\edge}[1]{\mbox{$#1$}}
\newcommand{\vertex}[1]{\mbox{\tiny #1}}
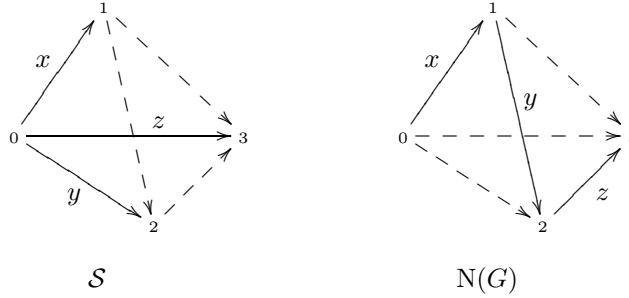
\[
\xymatrix@C=9pt@R=9pt{
& & \vertex{1} \ar@{-->}[dddddr] \ar@{-->}[dddrrr] & & & \\
& & & & & \\
& & & & & \\
\vertex{0} \ar[uuurr]^{\edge{x}} \ar[rrrrr]^(.63){\edge{z}} \ar[rrrdd]_{\edge{y}} & & & & & \vertex{3} \\
& & & & & \\
& & & \vertex{2} \ar@{-->}[uurr] & &
}
\hspace{5em}
\xymatrix@C=9pt@R=9pt{
& & \vertex{1} \ar[dddddr]^(.45){\edge{y}} \ar@{-->}[dddrrr] \\ \\ \\
\vertex{0} \ar[uuurr]^{\edge{x}} \ar@{-->}[rrrrr] \ar@{-->}[rrrdd] & & & & & \vertex{3} \\ \\
& & & \vertex{2} \ar[uurr]_{\edge{z}}
}
\]
\[
\xymatrix@C=9pt@R=9pt{& & \S & & &}
\hspace{5em}
\xymatrix@C=9pt@R=9pt{& & \mathrm{N}(G) & & &}
\]
Note what happens if we forget about the $0$:th vertex from the left tetrahedron, that is applying the d{\'e}calage:
The remaining vertices $(1, 2, 3)$, now resembling a 2-simplex, may be identified with the triple $(x, y, z)$ and their simplicial structure has nothing to do with the group structure. This describes the identification $$\Sigma \S = \NN_+(M)$$
for 2-simplices, that is
$(\Sigma \S)_\ord{2} = \S_\ord{3} = M^{\times 3} = \NN_+(M)_\ord{2}$.

\section{Conclusion}\label{conclusion}
We have shown that:
\begin{itemize}
 \item A simplicial triangulation on a sesquicategory $\C$ is the same thing as a weak recalage of its augmented lax nerve $\LN_+(\C)$.
 \item A group structure on a set $M$ is the same thing as a recalage of its augmented nerve $\NN_+(M)$.
 \item Every row in the following diagram is the categorification of the one below and that the upper row is the double categorification of the lower:
\end{itemize}
\vspace{1em} 
\begin{center}
\begin{tabular}{ l || l | l }
   Underlying object & Axiomatic structure & Structure as representaion \\
   \hline
   Sesquicategory & Simplicial triangulation & Weak recalage of lax nerve \\
   Category & Semitriangulation & \\
   Set & Group structure & Recalage of nerve \\
\end{tabular}
\end{center} 
\vspace{1em} 

The empty place is empty because of the non-functoriallity of the mapping cone.
We claim that a weak recalage of $\LN_+(\C)$ is an improvement over corresponding semitriangulation on $\Ho(\C)$.
Firstly, the definition of a recalage is a one-liner compared to a whole page for semitriangulation.
Secondly, it is an enhancement in the sense that Axiom {\bf S2} \& {\bf S4} become a functor (the mapping cone) at the cost of reformulating, rather than loosing the other axioms.

Moreover, in the case $\C = \Com{2}(\A)$ we have:
\begin{itemize}
 \item The weak recalage of $\LN_+(\Com{2}(\A))$ can be replaced with a strict one, namely the simplicial category $\Com{\Delta}(\A)$, with $\Com{\Delta}(\A)_\ord{1} = \Com{2}(\A)$.
\end{itemize}

This last point is perhaps the most interesting. It says that the asymmetric (the irregular definition of the $-1$ indices) and weak (simplicial relations only up to equivalence) properties of the recalage of $\LN_+(\Com{2}(\A))$ given by Example \ref{sim_tr_ch_com} can be overcome by the simplicial category $\Com{\Delta}(\A)$.
The simplicial maps $\mathrm{\bf d}_i$ and $\mathrm{\bf s}_i$ of $\Com{\Delta}(\A)$ behave in a symmetric way and the simplicial relations are satisfied up to equality rather than equivalence.
Moreover, there is a simplicial functor between their d{\'e}calage
$$\Spine : \Sigma \Com{\Delta}(\A) \longrightarrow \LN_+(\Com{2}(\A))$$
which is a pointwise adjoint equivalence.
But it goes only in one direction, the collection of pointwise weak inverses $\Filler_\ord{i}$ do not constitute a simplicial functor in the other direction.
This suggests that something important is being lost by $\Spine$ and that there is a better candidate for the categorification of semitriangulation of the homotopy category of chain complexes $\Ho(\Com{2}(\A))$, namely the simplicial category of $N$-complexes $\Com{\Delta}(\A)$.

\section*{Acknowledgements}

I am heartily thankful to my advisor, professor Volodymyr Mazorchuk.


\begin{thebibliography}{xxxxx}

\bibitem[Ab]{1}
\textsc{V. Abramov} On a graded $q$-differential algebra. Journal of Nonlinear Mathematical Physics 13, Suppl. 1, 1--8, 2006.

\bibitem[Bi]{2}
\textsc{J. Bichon} $N$-complexes et alg{\`e}bres de Hopf. Comptes Rendus Mathematique, Vol. 337, Issue 7, 2003.

\bibitem[CSW]{3}
\textsc{C. Cibils, A. Solotar, R. Wisbauer} $N$-complexes as functors, amplitude cohomology and fusion rules. Communications in Mathematical Physics, Vol. 272, Issue 3, 2007.

\bibitem[DV1]{4}
\textsc{M. Dubois-Violette} Tensor Product of $N$-complexes and Generalization of Graded Differential Algebras. Bulg. J. Phys. 36, No. 3, 227-236, 2009.

\bibitem[DV2]{5}
\textsc{M. Dubois-Violette} $d^N = 0$. arXiv:q-alg/9710021, 1997.

\bibitem[DV3]{6}
\textsc{M. Dubois-Violette} Generalized differential spaces with $d^N = 0$ and the $q$-differential calculus. Czechoslovak Journal of Physics, Vol. 46, Issue 12, 1996.

\bibitem[DVH]{7}
\textsc{M. Dubois-Violette, M. Henneaux} Tensor fields of mixed Young symmetry type and $N$-complexes. 
Communications in Mathematical Physics, Vol. 226, Issue 2, 2002.

\bibitem[DVK]{8}
\textsc{M. Dubois-Violette, R. Kerner} Universal $q$-differential calculus and $q$-analog of homological algebra. Acta Mathematica Universitatis Comenianae, Vol. 65, No. 2, 1996.

\bibitem[Es]{9}
\textsc{S. Estrada} Monomial Algebras over Infinite Quivers. Applications to $N$-Complexes of Modules. Communications in Algebra, Vol. 35, Issue 10, 2007.

\bibitem[GJ]{GJ}
\textsc{P. G. Goerss, J. F. Jardine.} Simplicial Homotopy Theory. Birkh\"{a}user, 1999.

\bibitem[GM]{GM}
\textsc{S. I. Gelfand, Y. I. Manin.} Methods of Homological Algebra. Springer, 2003.

\bibitem[Ha]{Ha}
\textsc{M. Hall, Jr.} The Theory of Groups. Chelsea, 1976.

\bibitem[He]{10}
\textsc{M. Henneaux} $N$-Complexes and Higher Spin Gauge Fields. International Journal of Geometric Methods in Modern Physics, Vol. 5, Issue 8, 2008.

\bibitem[IKM]{IKM}
\textsc{O. Iyama, K. Kato, J. Miyachi.} Derived categories of $N$-complexes. arXiv:1309.6039, 2014

\bibitem[Ka]{Ka}
\textsc{M. M. Kapranov.} On the $q$\textendash analog of homological algebra. arXiv:q-alg/9611005, 1991.

\bibitem[KW]{11}
\textsc{C. Kassel, M. Wambst} Alg{\`e}bre homologique des $N$-complexes et homologie de Hochschild aux racines de l'unit{\'e}. Publ. Res. Inst. Math. Sci. 34, No. 2, 91--114, 1998.

\bibitem[KQ]{KQ}
\textsc{M. Khovanov, Y. Qi.} An approach to categorification of some small quantum groups. arXiv:1208.0616, 2013.

\bibitem[Le]{Le}
\textsc{T. Leinster.} Higher Operads, Higher Categories. Cambridge University Press, 2004.

\bibitem[Mac]{Mac}
\textsc{S. Mac Lane.} Categories for the Working Mathematician, Second Edition. Springer, 1998.

\bibitem[May]{May}
\textsc{J. P. May.} The Axioms for Triangulated Categories. University of Chicago, Http, 2005.

\bibitem[Maye]{Maye}
\textsc{W. Mayer.} A new homology theory I, II. Annals of Mathematics, Vol. 43 p. 370--380, 594--605, 1942.

\bibitem[Maz]{Maz}
\textsc{V. Mazorchuk.} Lectures on algebraic categorification. arXiv:1011.0144, 2010.

\bibitem[Mi]{Mi}
\textsc{Dj. Mirmohades.} $N$--Complexes. Uppsala Universitet, diva-125738, 2010.

\bibitem[Ta]{12}
\textsc{A. Tikaradze} Homological constructions on $N$-complexes. Journal of Pure and Applied Algebra, Vol. 176, Issues 2–3, 2002.

\bibitem[Wal]{Wal}
\textsc{C. T. C. Wall.} On the Exactness of Interlocking Sequences. L'Enseignement Math{\'e}ma-tique, Band 12, 1966.

\bibitem[Wam]{13}
\textsc{M. Wambst} Homologie cyclique et homologie simpliciale aux racines de l'unit{\'e}. K-Theory, Vol. 23, 2001.

\bibitem[We]{We}
\textsc{C. A. Weibel.} An introduction to homological algebra. Cambridge University Press, 1994.

\bibitem[Wi]{Wi}
\textsc{N. H. Williams.} On Grothendieck universes. Compositio Mathematica, Vol. 21, Issue 1, 1969.

\bibitem[$n$Lab]{NL}
\textsc{$n$Lab} \href{http://ncatlab.org}{http://ncatlab.org}.

\end{thebibliography}
\end{document}